\documentclass[final,times,10pt]{article}
\usepackage{enumerate}
\usepackage{epsfig,alltt}
\usepackage[english]{babel}
\usepackage{blindtext}
\usepackage{dsfont}
\usepackage{amssymb}
\usepackage{amsfonts}
\usepackage{graphicx}
\usepackage{amsmath}
\usepackage{amsthm}
\usepackage[colorlinks=true,linkcolor=blue,citecolor=blue,pdfborder={0 0 0}]{hyperref}
\usepackage{natbib}
\usepackage[utf8]{inputenc}
\usepackage{lineno}
\usepackage[margin=1in]{geometry}

\font\dsrom=dsrom10 scaled 1200
\def \ind{\textrm{\dsrom{1}}}
\usepackage{soul}
\usepackage{floatflt}
\usepackage{makeidx}
\usepackage{algorithm}

\usepackage{algorithmic}

\usepackage{url}
\usepackage{float} 
\theoremstyle{remark}

\usepackage{xcolor}
\usepackage[tikz]{bclogo}

\usepackage{color}

\newcommand{\begeq}[1]{\begin{equation} \label{#1}}
	\newcommand{\fineq}{\end{equation}}

\newtheorem{prop}{\bf Proposition}
\newtheorem{lem}{ \bf Lemma}
\newtheorem{coro}{\bf Corollary}

\newtheorem{rem}{$\blacktriangleright$ \textbf{Remark}}

\title{Dependent censoring with simultaneous death times based on the Generalized Marshall-Olkin model}

\author{Mikael ESCOBAR-BACH$^{1}$\footnote{e-mail adress: mikael.escobar-bach@univ-angers.fr} and Salima HELALI$^{1}$\footnote{e-mail adress: helali.salima@gmail.com} \\
	$^{1}$Univ Angers, CNRS, LAREMA, SFR MATHSTIC, F-49000 Angers, France}

\numberwithin{equation}{section}
\begin{document}

    \baselineskip=25pt
	\date{ }
	\maketitle

\begin{abstract} 
In this paper, we considered the problem of dependent censoring models with a positive probability that the times of failure are equal. In this context, we proposed to consider the Marshall-Olkin type model and studied some properties of the associated survival copula in its application to censored data. We also introduced  estimators for the marginal distributions and the joint survival probabilities under different schemes and showed their asymptotic normality under appropriate conditions. Finally, we evaluated the finite-sample performance of our approach relying on  a small simulation study on synthetic data, and  an application to real data.\\
\end{abstract}

\section{Introduction}
The development of new approaches for partially recorded data has become a major concern in several statistical areas due to right-censoring issues that frequently arise in statistical applications. The phenomenon occurs when the follow-up time of interest $T$ is eventually shortened by a random censoring time $C$ and where only $Y=\min(T,C)$ is given with the information that $Y=T$ or $Y=C$. In the literature, most papers require the assumption that the survival time $T$ is independent of the censoring time $C$, which is rarely met in practice. For instance, patients suffering from a vascular disease might die from either a heart attack or a stroke, which are likely to be correlated and censored from one if the other arises. Although the dependent censoring model is thus more realistic, it has been shown  that it lacks the nonparametric identifiability for the joint distribution of $(T,C)$ in the absence of prior knowledge of its dependence structure, \cite{Tsia75}. In this regard, \cite{Zhe95} introduced a copula-graphic estimator for the marginal distributions when the copula is known. Later, \cite{Riv01} further investigated this proposal by assuming an Archimedean copula. This approach has been extended for various regression models imposing a known copula/archimedean generator function in \cite{Bra05}, \cite{Chen10}, \cite{de17}, \cite{Suj18} and \cite{Emu18} among others.  Some other recent developments managed to relax the assumption of known copula by considering parametric models in \cite{Cza22} and semiparametric models in \cite{Der23a}.\\
Similar to the censoring framework, competing risk models arise when individuals are at risk for more than one mutually exclusive event. The Marshall-Olkin (MO) shock model introduced in \cite{Mar67} is a competing risk model with an alternative expression concerning the follow-up times. Here, the lifetimes are subject to three independent types of shocks $(X_1,X_2,X_3)$ such that $T=\min(X_1,X_3)$ and $C=\min(X_2,X_3)$. The model was originally defined to handle complex risk systems with exponential distributions where the two components operate in a random environment and are subject to a third common termination source. In the literature, several generalizations of MO-type distributions have been developed where the Weibull distribution is considered for the random shocks, see e.g. \cite{Mul87} and several recent contributions by \cite{Kun10,Kun13}, \cite{Mir15}, \cite{Moh14} and \cite{Sho12}, where parametric estimators have been proposed under maximum likelihood or Estimation-Maximization approaches. Further works have also relaxed the parametric margins assumption as in \cite{Li11} with a generalized Marshall–Olkin (GMO) model, where the random variables $X_i$ have an arbitrary continuous distribution. In a survival analysis extension, \cite{Li12} consider a Weibull MO-type model and propose two estimations procedures for bivariate censored data. For more insight on the classical MO model, we refer the reader to a recent review of the subject in \cite{Bay22}.\\ 
When it comes to consider credit risk modeling, competing events are likely to be equal when exposed to similar risks. As an example, a life insurance policy is a bounding contract between the  insured and the insurance company which protects the holders from their risks of death. If the contract is taken out by a married couple, the accrued capital is due to return to the remaining member when the other one dies. However, the partners are exposed to the same risks, such as catastrophic exogenous events, which may lead to the simultaneous death of the spouses. The dependency between these lifetimes has been shown in \cite{Carr00} and dependent random time models have been largely developed in the literature of credit risk and life insurance. The actuarial community usually agrees on two approaches; one with a Markovian representation, which handles the so-called {\it broken-heart} syndrome, referring to a change point in the lifetimes \citep[see e.g.][]{Den01,Spr13,Dic20} and another one with copula based models which allow to model the joint mortality \citep[see e.g.][]{Carr00,Spr06}. A recent work in \cite{Gob21} attempts to combine these two approaches and propose a Ryu-type model based on a MO probabilistic construction.\\ 

The work here focused on the analysis of censoring models where the probability of the joint default between the random lifetimes is strictly positive. Such line of work is important for the applicability of the survival models, in particular with dependent censoring. This, however, goes beyond the scope of copula based models with standard statistical assumptions and requires copula functions with singular components. We can cite several works studying the singular components of copulas, as in \cite{Dur15} and \cite{Fer15}, or proposing new copula constructions, as in \cite{Xie19} where the authors proposed a model unifying the Marshall-Olkin and Archimedean copulas. In this study, we  considered a dependent censoring scenario where $(T,C)$ is based on the GMO model, which has the sought-after property that $\mathbb{P}(T=C)>0$ when $X_3$ is non-degenerate. Such probabilistic assumption particularly induces dependent random lifetimes but makes difficult the dependence structure identification. Likewise \cite{Pro22}, we also extended the recent works in risk modelling linked to the MO shock models, by considering censored outputs where only the first occurring event time is recorded. With this regard, we first provided the model specifications, linked to the aforementioned works, and discussed the model identification depending on the information  of the censoring status. We  showed that we can adapt common statistical methods from survival analysis for estimating the different model specifications, and more particularly, the marginal distributions and/or the joint distribution.\\
The remainder of the paper is organized as follows.  In Section \ref{sec:model}, we identified  the model and the notations. We  introduced the main characteristics of the survival copula and its asymptotic properties.  In section \ref{sec:estimation},  we  suggested  non parametric estimators of  the survival distributions, the survival joint distribution and the Kendall's tau measure and established theirs asymptotic normality. Section \ref{sec:num} was devoted to the application results obtained, firstly through simulations (subsection \ref{sub1}) and secondly using a real data set (subsection \ref{sub2}). Some pertinent  concluding remarks and potential perspectives  were provided in Section \ref{sec:con},  while Section \ref{sec:pro}  offered the proofs of the theoretical results.

\section{Model and survival copula}\label{sec:model}
The GMO model is a direct extension of the exponential MO model for random variables with arbitrary distributions. It is defined as the competing risk between three independent and positive random variables $X_1$, $X_2$ and $X_3$. We denoted their respective distributions and cumulative hazard functions by $F_i$ and $\Lambda_i$ for $i=1,2,3$. In our context,  it is worth recalling that the lifetimes of interest are
\begin{eqnarray}
\label{gmomodel}
    T=\min(X_1,X_3)\quad\text{and}\quad C=\min(X_2,X_3).
\end{eqnarray}
It is clear that $T$ and $C$ are dependent as they both share the same failure time $X_3$. Due to the censoring mechanism, we considered that the only observed time is given by $Y=\min(T,C)$,  whose distribution function is denoted by $H$. Recall that for any generic distribution function $F$, its cumulative hazard function is given by
\begin{eqnarray*}
\Lambda(t)=\int_0^t\dfrac{F(du)}{\overline{F}(u^-)},\quad t\geq 0,
\end{eqnarray*}
where $\overline{F}:=1-F$ and $\overline{F}(u^-)=\lim_{s\downarrow u}\overline{F}(s)$. Reversely, the relationship between a distribution function $F$ and its cumulative hazard function $\Lambda$ is given by
\begin{eqnarray}
\label{eq::lambda}
\overline{F}(t)=:\exp\left(-\widetilde{\Lambda}(t)\right)=\exp\left(-\Lambda^c(t)+\sum_{s\leq t}\log(1-\Delta\Lambda(s))\right),    
\end{eqnarray}
where $\Lambda^c$ and $\Delta\Lambda$ respectively denote the continuous and discontinuous parts of $\Lambda$ \citep[see p.898 in][]{Sho86}. As such, the model can be described through the survival functions 
\begin{eqnarray*}\label{marginal}
&&\overline{F}_T(t)=:\exp(-\widetilde{\Lambda}_T(t))=\exp(-\widetilde{\Lambda}_1(t)-\widetilde{\Lambda}_3(t))\\
\text{and}&&\overline{F}_C(t)=:\exp(-\widetilde{\Lambda}_C(t))=\exp(-\widetilde{\Lambda}_2(t)-\widetilde{\Lambda}_3(t)),
\end{eqnarray*}
where $\Lambda_T$ and $\Lambda_C$ define the cumulative hazard functions of $T$ and $C$ respectively. Similarly, the joint survival function follows with
\begin{eqnarray*}
    \widetilde{P}(t,s)=\mathbb{P}(T>t,C>s)=\exp(\widetilde{\Lambda}_1(t)-\widetilde{\Lambda}_2(t)-\widetilde{\Lambda}_3(\max(t,s)))\quad t, s \geq 0.
\end{eqnarray*}
Note that bivariate distributions with this latter form characterize the GMO model as introduced in \cite{Li11} where the functions $\widetilde{\Lambda}$ or called \textit{generating functions} in their work. We next proposed to describe the dependence structure of the random vector $(T,C)$ thanks to their copula function and  examine some of their properties. A copula stands for a bivariate distribution function $\mathcal{C}$ with uniform margins. \cite{Skl59} elucidated the role that copulas play in the inner relationship for multivariate distributions and showed that they enable  free-margins description for dependence structures. We know that for any marginal distribution functions $F_T$ and $F_C$,  there exists  only one  copula $\mathcal{C}$,  determined on the set $F_T(\mathbb{R}_+)\times F_C(\mathbb{R}_+)$, for which
\begin{eqnarray*}\label{copule}
P(t,s)=\mathbb{P}(T\leq t,C\leq s)= \mathcal{C}\left(F_T(t),F_C(s)\right),\quad  t, s \geq 0.
\end{eqnarray*}
In reliability studies, although this would not change the analysis, it is often more convenient to express a joint survival function in terms of its marginal survival functions, such that
\begin{eqnarray}\label{copulesurvie}
\widetilde{P}(t,s)=\widetilde{\mathcal{C}} \left(\overline{F}_T(t),\overline{F}_C(s)\right),\quad t, s \geq 0,
\end{eqnarray}
where the function $\widetilde{\mathcal{C}}$ is called the \textit{survival copula}. In \cite{Li11}, a general form for the survival copula is derived and used to study some aging properties. We here proposed another form for the dependence structure, intended to generalize the form of the exponential MO model.

\begin{prop}\label{prop::generalform}
The survival copula for the GMO model is given by
\begin{eqnarray}\label{survivalcopula}
\widetilde{\mathcal{C}}(u,v)=uv\min\left(u^{-\alpha_1\left(\overline{F}^{-}_T(u) \right)},v^{-\alpha_2\left(\overline{F}^{-}_C(v) \right)}\right),\quad
(u,v)\in \overline{F}_T(\mathbb{R}_+)\times \overline{F}_C(\mathbb{R}_+),
\end{eqnarray}
where $\displaystyle{\alpha_i=\widetilde{\Lambda}_{3}/(\widetilde{\Lambda}_i+\widetilde{\Lambda}_{3}})$ for $i=1,2$ and $\overline{F}^{-}$ stands for the generalized inverse function of $\overline{F}$.\\
\end{prop}

\begin{rem} In the sequel, we will say that a MO survival copula is  defined as in (\ref{survivalcopula}) with constant functions $\alpha_i$ for $i=1,2$, which we sometimes call parameters for the ease of reading. In this case the Kendall's tau  is given in terms of $\alpha_i$ by
$\displaystyle{\tau=\frac{\alpha_1 \alpha_2}{\alpha_1-\alpha_1\alpha_2+\alpha_2}}.$\\
\end{rem}

Note that the minimum part in (\ref{survivalcopula})  reflects the discrepancy for the copula away from the independent scenario, obtained when $\alpha_i\equiv 0$,  however the complete dependence case is fulfilled when $\alpha_i\equiv 1$. In a more general  way, constant functions $\alpha_i$ are directly linked to the MO model. As a consequence  of Proposition \ref{prop::generalform}, it turns out that the GMO model  fits for proportional shock events only.
\begin{coro}
\label{cor::prop}
    Let $(T,C)$ be any bivariate distribution defined as in (\ref{gmomodel}). Then, its survival copula function is given by (\ref{survivalcopula}) with constant functions $\alpha_i$ for $i=1,2,3$ if and only if the functions $\{\widetilde{\Lambda}_i,\,i=1,2,3\}$ are proportional.
\end{coro}
It is sometimes more convenient to summarize the model complexity by  a universal measure of the dependency, such as the Kendall's tau correlation. In our next procedure, the Kendall’s tau corresponding to the GMO model was based on \cite{Mul18} finding, but with an alternative expression which would be useful for the measure estimation.

\begin{prop}\label{prop::kendall}
    Let us assume that $X_i,\,i=1,2,3$, have cumulative hazard functions that are continuous and strictly increasing. We also assume that $\lim\limits_{t\to +\infty}\max\{\Lambda_i(t),\,i=1,2,3\}=+\infty$. Then, the Kendall's tau associated to the GMO model is given by
\begin{eqnarray}\label{eq::tau}
\tau=2\int_0^{+\infty}\overline{H}(u)H_3^1(du),
\end{eqnarray}
where $H_3^1$ denotes the sub-distribution function given by $H_3^1(t):=\mathbb{P}(Y\leq t, T=C)$, $t>0$.
\end{prop}
From a copula-based model perspective, it is also interesting to determine the conditions ensuring that a bivariate distribution with an MO survival copula admits an MO probabilistic construction as in (\ref{gmomodel}). In the next result, we showed that a generalization of this property  would be possible under a proportional hazard condition.  

\begin{lem}
\label{lem::gmo}
Let $(X,Y)$ be a random vector drawn from  an MO survival copula with parameters $\alpha_1$ and $\alpha_2$ in $(0,1]$. Then, there exists random independent variables $\{X_i,\,i=1,2,3\}$ such that the random couple 
\begin{eqnarray*}
    (T,C)=\left(\min(X_1, X_3),\,\min\left(X_2,\widetilde{\Lambda}^-_Y\left(\frac{\alpha_1}{\alpha_2}\widetilde{\Lambda}_X(X_3)\right)\right)\right),
\end{eqnarray*}
follows the same law as $(X, Y)$. In particular, the distribution of $(X,Y)$ can be interpreted as a GMO model if $\widetilde{\Lambda}_Y=\frac{\alpha_1}{\alpha_2}\widetilde{\Lambda}_X$.
\end{lem}
 It is also revealed in  \cite{Li11} (see Theorem 4.5) that  $(T,C)$  and its  residual life $(T_t,C_t):=(T-t,C-t| T>t,C>t)$  share the same GMO copula  if and only if $\widetilde{\Lambda}_i$ for $i=1,2,3$ are proportional. This lack-of-memory property particularly suggests that the survival copula for the GMO model should belong to a maximum domain of attraction of some extreme value copula.  It should be recalled here that a copula $C$ is defined as an extreme value copula if there exists an alternative copula $C_I$ such that
\begin{eqnarray*}
\lim_{n\to+\infty}C_I\left(u_1^{1/n},u_2^{1/n}\right)^{n}=C(u_1,u_2),
\end{eqnarray*}
and the copula $C_I$ is said to belong to maximum domain of attraction of $C$. It is reversely clear that a survival copula in (\ref{survivalcopula}) with constant functions $\alpha_i$ are extreme value copula. To show that the GMO model is in a maximum domain of attraction, we considered the distributions $F$ under the \textit{von Mises'} condition \citep[see e.g.][]{De06}, that is   
\begin{eqnarray}\label{cond1}
\lim\limits_{t \to x^{\star}_F} \left( \frac{1-F}{F'} \right)' (t)=\gamma,
\end{eqnarray}
for some $\gamma\in\mathbb{R}$ and where $x^{\star}_F$ defines the $F$ terminal point. The distribution is then said to belong to the maximum domain of attraction $\mathcal{D}_\gamma$.
\begin{prop}\label{carac2}
Suppose that the $F_i$'s satisfies the \textit{von Mises'} condition for some $\gamma_i$, all of the same sign but with possibly different values. Let $\xi_i=\gamma_i/(\gamma_i+\gamma_{3})$ for any $i=1,2$ and define the extreme value copula $C_e(u,v):=\min(vu^{1-\xi_1},uv^{1-\xi_2})$ with Pickands dependence function given by $\displaystyle{A(t):=1-\min(\xi_1(1-t),\xi_2 t)}$. Then :
\begin{enumerate}[i)]
\item if $\gamma_i >0$, then  $x^{\star}_{F_T}=x^{\star}_{F_C}=\infty$ and $\widetilde{C}$ belongs to the maximum domain of attraction of $C_e$.
\item if $\gamma_i<0$ and $x^{\star}_{F_1}=x^{\star}_{F_2}=x^{\star}_{F_{3}}$, then 
$x^{\star}_{F_T}=x^{\star}_{F_C}<\infty$ and $\widetilde{C}$ belongs to the maximum domain of attraction of $C_e$.
\end{enumerate}
\end{prop}
\noindent
Note that the lighted-tail case with $\gamma_i=0$ can also gives extreme value copula but requires second order assumptions on the tail distributions which greatly complicates the analysis and would hence be consider in further works.

\section{Estimation}\label{sec:estimation}

\subsection{Cumulative Hazard functions and marginal distributions}
Unlike most survival models with dependent censoring, it turns out that the Nelson-Aalen estimator  $\Lambda_{T,n}$ of the cumulative hazard function $\Lambda_T$   remains consistent based on the considered model of dependent censoring (\ref{gmomodel}).
In a more general manner, on can show that the Nelson-Aalen estimator can reach for the different cumulative hazard functions $\Lambda_T$, $\Lambda_C$ or $\Lambda_i$ for $i=1,2,3$ which characterize the GMO if the expression of the event of interest  $Y=T$, $Y=C$ or $Y=X_i$ is observed respectively.  The following proposition shows how to rewrite the cumulative hazard functions according to the sub-distribution of the associated event.
\begin{prop}\label{Lambda1etT}
Let $\delta^{(k)}=\ind_{\cap_{j\neq k}\{ X_k \leq X_j\}}$ with $k=1,2,3$. Then if we define the sub-distribution functions $H_k^1(t):=\mathbb{P}(Y\leq t,\delta^{(k)}=1)$, we have
\begin{eqnarray*}
\Lambda_k(t)= \int_0^t \frac{dH^1_k(u)}{ \overline{H}(u^-)},\quad t\geq 0.  
\end{eqnarray*}
Moreover, if $\mathbb{P}(X_1=X_3)=0$ or $\mathbb{P}(X_2=X_3)=0$, we respectively have 
\begin{eqnarray*}
\Lambda_T(t)= \int_0^t \frac{dH^1_4(u)}{\overline{H}(u^-)}\quad\text{and}\quad\Lambda_C(t)= \int_0^t \frac{dH^1_5(u)}{\overline{H}(u^-)},
\end{eqnarray*}
where $\delta^{(4)}=\ind_{\{T\leq C\}}$ and $\delta^{(5)}=\ind_{\{T\geq C\}}$.
\end{prop}
The previous proposition particularly emphasizes that the hazard functions of $T$ and $C$ keep the same form as in the case of independent censoring. Based on this property, we obtain a nonparametric estimators of the cumulative hazard functions, which have the same expression as  the common Nelson-Aalen estimator, consistent even  under the dependent censoring model (\ref{gmomodel}). To do so, we define $\{(Y_i,\delta^{(k)}_i)\}_{1\leq i\leq n}$ as an i.i.d. sample drawn from the random couple $(Y,\delta^{(k)})$, where $k\in [\![1,5]\!]$ represents the known event. For any $i=1,\ldots,n$, also define $Y_{(i)}$ as the $i$-th order statistic  together with its associated indicator $\delta^{(k)}_{(i)}$. Based on (\ref{eq::lambda}), the Nelson-Aalen type estimators are then defined by 
\begin{eqnarray}\label{estLambda}
\Lambda_{k,n}(t):=\sum_{i: Y_{(i)} \leq t} \frac{d_{k,i}}{n-i+1}\quad\text{with}\quad d_{k,i}:=\sum_{j=1}^n\ind_{\{Y_j=Y_{(i)}, \delta^{(k)}_j=1\}}\quad\text{and}\quad t>0.
\end{eqnarray}
Since  the approximation of the Kaplan-Meier (KM) estimator, defined as 
\begin{eqnarray*}
\overline{F}_{k,n}(t)=
\begin{cases}
\prod\limits_{i:Y_{(i)} \leq t} \left( 1- \frac{1}{n-i+1}\right)^{\delta_{(i)}^{(k)}} & \text{ if } t < Y_{(n)},\\
0 & \text{ if } t \geq Y_{(n)},
\end{cases}
\end{eqnarray*}
is obtained through the Nelson-Aalen estimator according to Lemma 1 of \cite{Bre74}, we conclude that the KM estimator also provides a consistent
estimator of the survival function in the case of dependent censoring (\ref{gmomodel}).  Note that the asymptotic properties for $\overline{F}_{k,n}$ are then straightforward to obtain from common results on the KM estimator under a survival model with independent censoring \citep[see e.g.][among others]{And12,Sho86,Mal92}.

\begin{rem}
    It is worth noting here that the estimators in (\ref{estLambda}) are defined for each underlying distribution $k\in [\![1,5]\!]$ independently. In the sequel, we will see that the joint law estimation requires at least the observation of two different events $\delta^{(k)}$ and $\delta^{(k')}$ to be effective. Under the model $(Y,\delta^{(4)})$, this particularly implies that we obtain a consistent estimator for $F_T$, although the dependent structure is not always identifiable. 
\end{rem}

The next proposition ensures the asymptotic normality for the  previous estimators as a random vector in terms of  event types and  time values.
\begin{prop}\label{asympnormLambda}
Let $\{t_i\}_{1\leq i\leq m}$ for $m\in\mathbb{N}^*$ a series of positive reals and assume that the cdf's $F_i$ are continuous. Then the random vector
\begin{eqnarray*}
   \sqrt{n} (\Lambda_{k,n}(t_i)-\Lambda_k(t_i);\,k\in[\![1,5]\!],\,i\in[\![1,m]\!]),
\end{eqnarray*}
is asymptotically normal towards a centered random normal vector $N:=(N_k(t_i);\,k\in[\![1,5]\!],\,i\in[\![1,m]\!])$ with covariance matrix $\Phi:=\left(\sigma_{k,l}(t_i,t_j);\,i,j\in[\![1,m]\!],\,k,k'\in[\![1,5]\!]\right)$ where for $k=l$
\begin{eqnarray*}
\sigma_{k,k}(t_i,t_j)=\int_0^{t_i\wedge t_j} \dfrac{d\Lambda_k(u)}{\overline{H}(u)},
\end{eqnarray*}
and for $k\neq l$
\begin{eqnarray*}
\sigma_{k,l}(t_i,t_j)&=&\int_0^{t_j} \int_0^{t_i}  \frac{\overline{H}(u\vee v)}{\overline{H}(u)\overline{H}(v)}d\Lambda_k(u) d\Lambda_l(v)- \int_0^{t_j \wedge t_i} \frac{\Lambda_k(t)-\Lambda_k(u)}{\overline{H}(u)}d\Lambda_l(u)
\\ & & - \int_0^{t_j\wedge t_i} \frac{\Lambda_l(s)-\Lambda_l(u)}{\overline{H}(u)}d\Lambda_k(u)+ \mathbb{E}\left[\frac{\ind_{\{ Y \leq t_j\wedge t_i, \delta^{(k)}=\delta^{(l)}=1\}}}{\overline{H}(Y)^2} \right].
\end{eqnarray*}
\noindent
Note that the last term can be null if we have mutually exclusive events between $\delta^{(k)}$ and $\delta^{(l)}$.
\end{prop}
\subsection{Joint distribution} \label{nonparasurvivalcopule}
In the previous section, we have seen that the estimation of the marginal distributions is possible as long as the correct event of interest is known. However, it has been shown that due to dependent censoring, the joint law distribution (\ref{copulesurvie}) is not always identifiable. In fact, a complete characterization of the GMO and its dependence structure remains impossible with only one indicator $\delta^{(k)}$ at hand,  i.e. the mapping between the distributions of the $X_i$'s and $(Y,\delta^{(k)})$ might be not one-to-one. For instance, let $X_i$ be an exponential random variable with parameter $\lambda_i$ for $i=1,2,3$ and consider its GMO model $(Y,\delta^{(4)})$. Then, for any $t\geq 0$
\begin{eqnarray*}
        &&\mathbb{P}(Y>t,\delta^{(4)}=1)=\frac{\lambda_1+\lambda_3}{\lambda_1+\lambda_2+\lambda_3}e^{-(\lambda_1+\lambda_2+\lambda_3)t},\\
        \text{ and }&&\mathbb{P}(Y>t,\delta^{(4)}=0)=\frac{\lambda_2}{\lambda_1+\lambda_2+\lambda_3}e^{-(\lambda_1+\lambda_2+\lambda_3)t},
\end{eqnarray*}
only depend on $\lambda_1+\lambda_3$ and $\lambda_2$, which are the parameters for the distributions of $T$ and $X_2$ respectively. Hence, the estimation of the survival copula and/or the joint distribution  requires the differentiation of at least two events. Within the survival analysis context, this information is available through the observation of the censoring indicators $(\delta^{(4)},\delta^{(5)})$ traducing whether the observed time $Y$ is censored $(Y=C)$ or not $(Y=T)$.  Thus, in order to overcome the problem of the identifiability  of the joint law distribution  and  be able to estimate it, we considered the i.i.d. sample $\left(Y_i,\delta^{(4)}_i,\delta_i^{(5)}\right)_{1\leq i\leq n}$ drawn from the model $\left(Y,\delta^{(4)},\delta^{(5)}\right)$. We then proposed an estimation for the GMO model based on (\ref{copulesurvie}) by replacing each function with their empirical counterparts, which  gave us an estimator for the survival bivariate function given by 
\begin{eqnarray*}
\widetilde{P}_n(t,s)=\min\left(\overline{F}_{5,n}(s)\overline{F}_{4,n}(t)^{1-\alpha_{1,n}\left(t \right)},\overline{F}_{4,n}(t)\overline{F}_{5,n}(s)^{1-\alpha_{2,n}\left(s \right)}\right),\quad \text{for }s,t\geq 0
\end{eqnarray*}
where
\begin{eqnarray*}
\alpha_{1,n}(t)=\frac{\Lambda_{3,n}(t)}{\Lambda_{4,n}(t)}\quad\text{and}\quad\alpha_{2,n}(s)=\frac{\Lambda_{3,n}(s)}{\Lambda_{5,n}(s)}. 
\end{eqnarray*}
In the next result, we finally derived the asymptotic properties for the joint distribution estimator.
\begin{prop}\label{asympnorm}
  Let $c>0$ such that $F_i(c)<1$ for $i=1,2,3$ and assume that the $F_i$'s are continuous. For any $0 \leq t,s \leq c$ such that $\overline{F}_T(t)^{\alpha_1(t)}\neq\overline{F}_C(s)^{\alpha_2(s)}$, we have 
  \begin{eqnarray*}
      \sqrt{n}\left(\dfrac{\widetilde{P}_n(t,s)}{\widetilde{P}(t,s)}-1\right) \underset{n \to +\infty}{\overset{\mathcal{L}}{\longrightarrow}}\mathcal{N}(0,\Sigma(t,s)) ,
  \end{eqnarray*}
  with
  \begin{eqnarray*}
    \Sigma(t,s)=
    \begin{cases}
       \sigma_{3,3}(s,s)+\sigma_{4,4}(t,t)+\sigma_{5,5}(s,s)&
       \\+\,2(\sigma_{4,5}(t,s)-\sigma_{4,3}(t,s)-\sigma_{5,3}(s,s)), & \text{if } \overline{F}_T(t)^{\alpha_1(t)}<\overline{F}_C(s)^{\alpha_2(s)}.\\[.2cm]
       \sigma_{3,3}(t,t)+\sigma_{4,4}(t,t)+\sigma_{5,5}(s,s)&\\
       +\,2(\sigma_{4,5}(t,s)-\sigma_{4,3}(t,t)-\sigma_{5,3}(s,t)), & \text{if }\overline{F}_T(t)^{\alpha_1(t)}>\overline{F}_C(s)^{\alpha_2(s)}.
    \end{cases}
  \end{eqnarray*}
    In the case $\overline{F}_T(t)^{\alpha_1(t)}=\overline{F}_C(s)^{\alpha_2(s)}$, 
    the weak convergence is true with
    \begin{eqnarray*}
           \sqrt{n}\left(\dfrac{\widetilde{P}_n(t,s)}{\widetilde{P}(t,s)}-1\right) \underset{n \to +\infty}{\overset{\mathcal{L}}{\longrightarrow}}\mathcal{Z}(t,s) ,
    \end{eqnarray*}
    where $\mathcal{Z}(t,s)$ is a random variable which follows the same law as the minimum between two random normal variables $N_3(t)-N_4(t)$ and $N_3(s)-N_5(s)$ where $N:=(N_3(t),N_4(t),N_3(s),N_5(s))$ is a random normal vector as given in Proposition \ref{asympnormLambda}.
 \end{prop}
 \subsection{Kendall's tau}
In this section, we proposed an estimator for the Kendall's tau measure based on  Proposition \ref{prop::kendall}. We considered the same sample from the previous section and constructed our estimator with the help of the representation (\ref{eq::tau}) by replacing the functions by their empirical versions. This gives us the following estimator 
\begin{eqnarray*}
    \tau_n=\dfrac{2}{n}\sum_{i=1}^n\delta_{(i)}^{(4)}\delta_{(i)}^{(5)}\dfrac{n-i+1}{n},
\end{eqnarray*}
where we recall that $\delta_{(i)}^{(4)}$ and $\delta_{(i)}^{(5)}$ denote the indicator associated to the $i$-th order statistic $Y_{(i)}$. It is worth noting here that $\delta^{(4)}\delta^{(5)}=\delta^{(3)}$ which indicates the simultaneous death times event ${T=C}$. In the next proposition, we showed the asymptotic normality for the Kendall's tau estimator under mild conditions.
\begin{prop}
\label{prop::kendallest}
Consider that the cdf's $F_i$ are continuous. If $x^*_{F_T}=x^*_{F_C}=\infty$, we assume that $H$ also satisfies the \textit{von Mises'} condition with extreme value index $0\leq \gamma<1/4$. We further consider that if $\gamma\neq0$, then there exists a bounded slowly varying function $\ell$ such that $H(t)=t^{-1/\gamma}\ell(t)$ for $t$ large enough. Then we have
$$\sqrt{n}(\tau_n-\tau) \underset{n \to +\infty}{\overset{\mathcal{L}}{\longrightarrow}}\mathcal{N}(0,\sigma^2),$$
where 
\begin{eqnarray*}
\sigma^2 &=&  4\int_0^{+\infty}\int_0^{+\infty} \overline{H}(x\lor z) - \overline{H}(x) \overline{H}(z) dH_3^1(x) dH_3^1(z)+ 4\int_0^{+\infty}\int_0^{+\infty} H_3^1(x \wedge z) - H_3^1(x) H_3^1(z) dH(x) dH(z)
\\& & + 8 \int_0^{+\infty} \int_0^{+ \infty} \mathbb{P} \left(x < Y \leq z,\delta^{(3)}=1 \right)-\overline{H}(x) H_3^1(z)  dH_3^1(x) dH(z).
\end{eqnarray*}
\end{prop}

\section{ Simulation and real data analysis}\label{sec:num}

\subsection{ A small simulation study}\label{sub1}
\begin{figure}[H]
\centering
\includegraphics[width =6cm, angle=0]{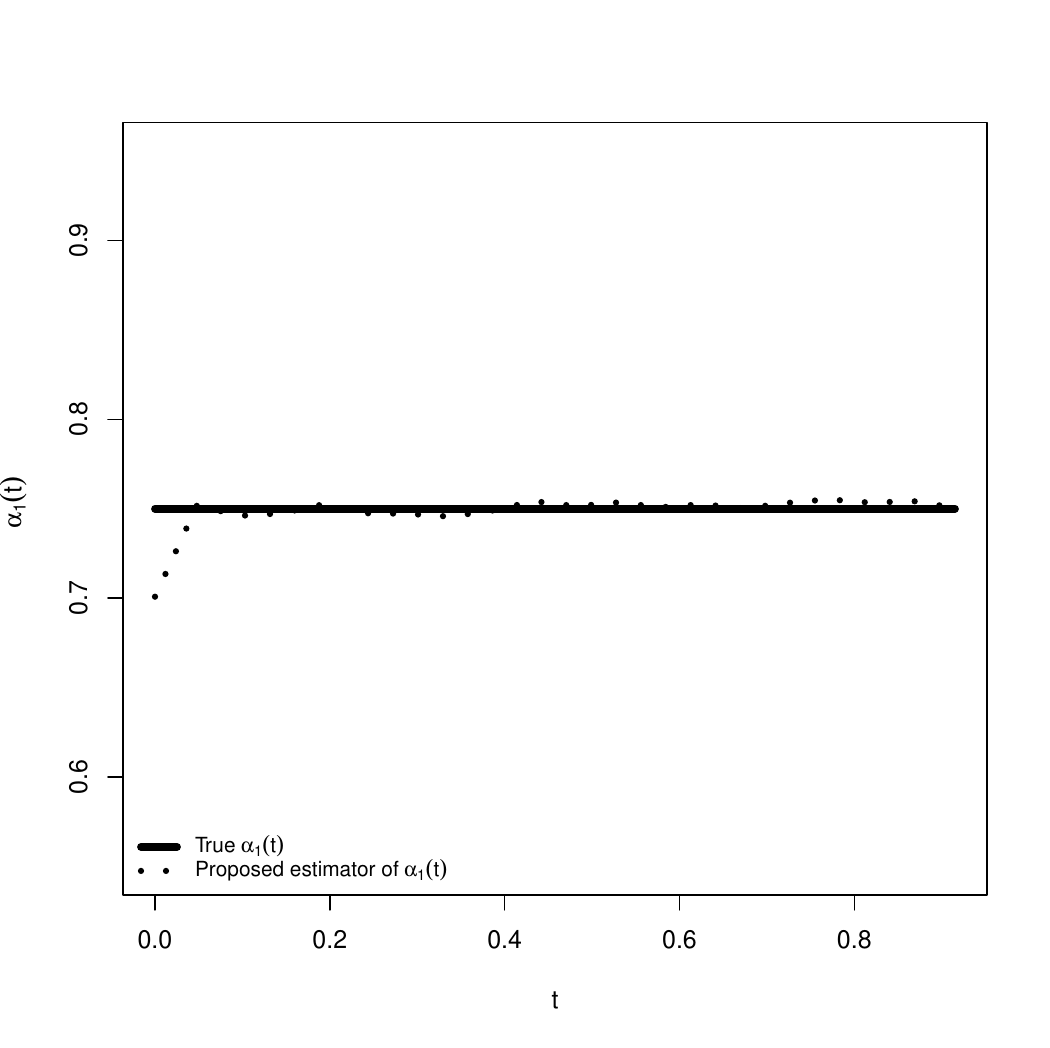}
\includegraphics[width =6cm, angle=0]{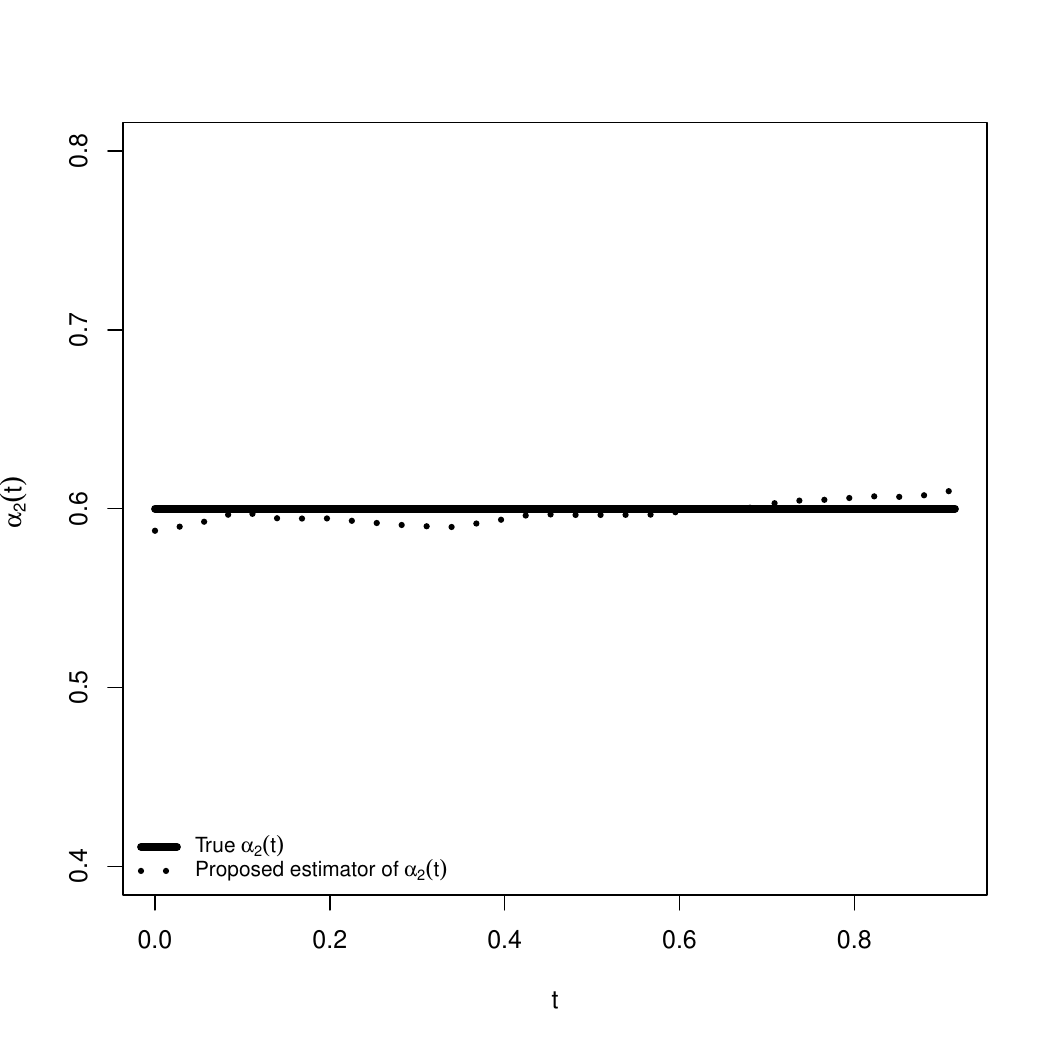}
\includegraphics[width =6cm, angle=0]{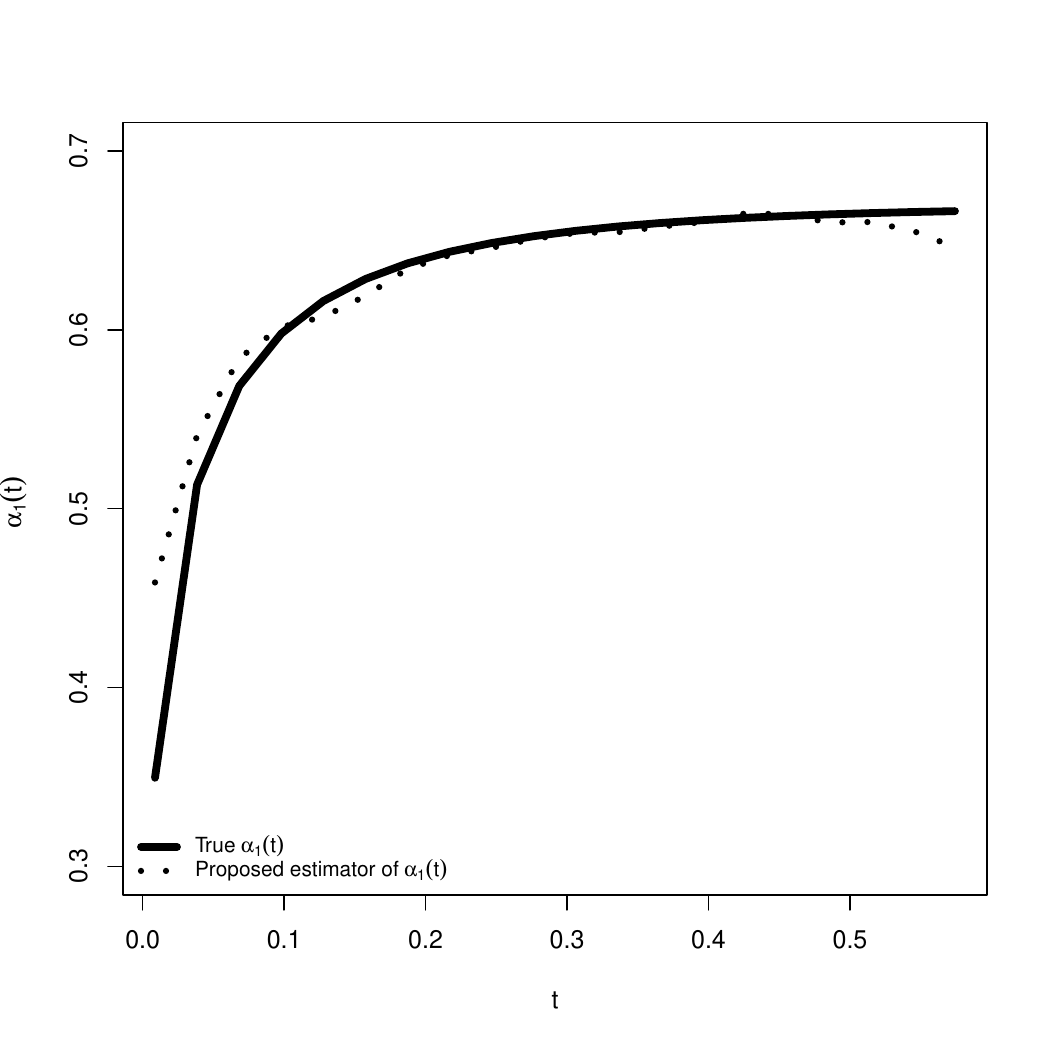}
\includegraphics[width =6cm, angle=0]{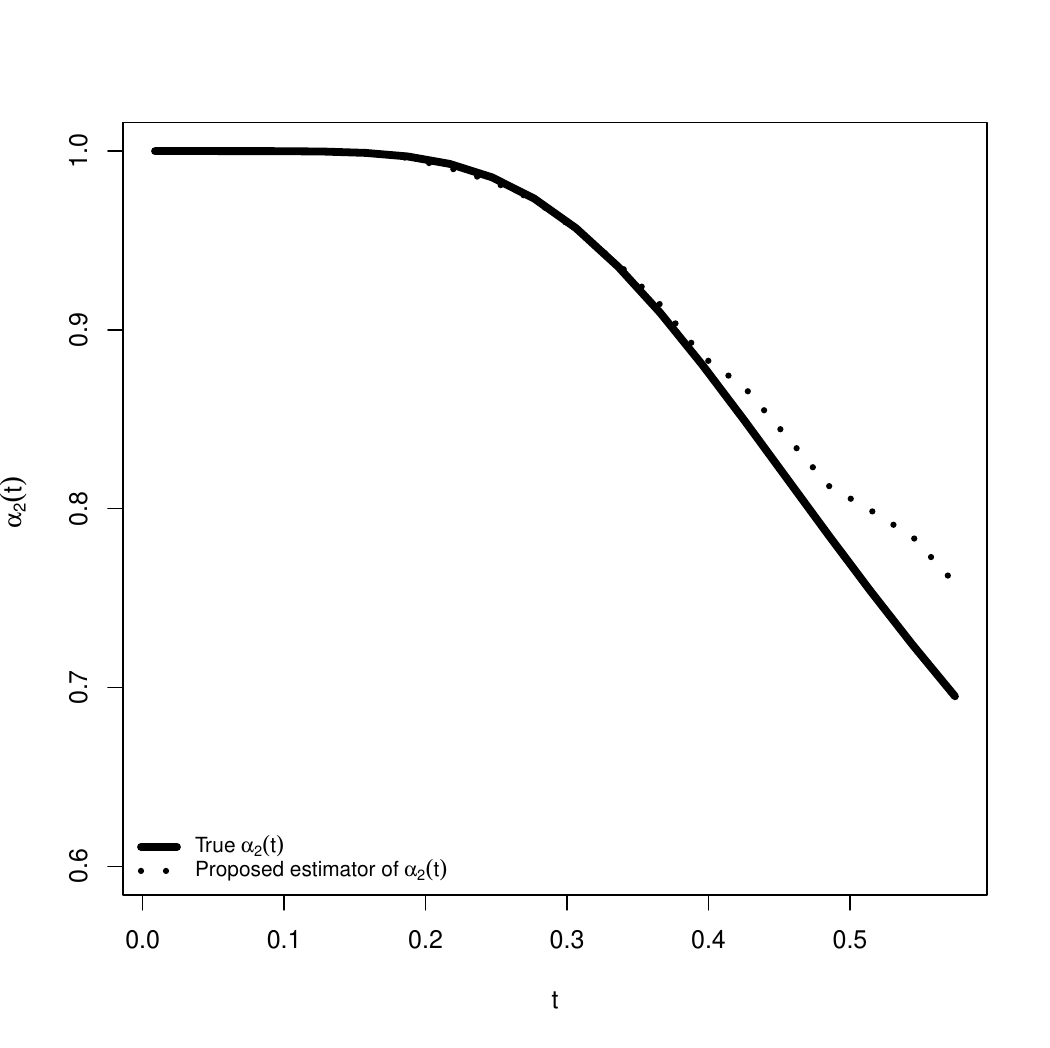}
\caption{Estimators of $\alpha_1$ (left) and $\alpha_2$ (right) for $N=100$ samples of size $n=200$ for the  model $a$ (up) and $b$ (down). }\label{alpha}
\end{figure}

\begin{figure}[H]
\centering
\includegraphics[width =6cm, angle=0]{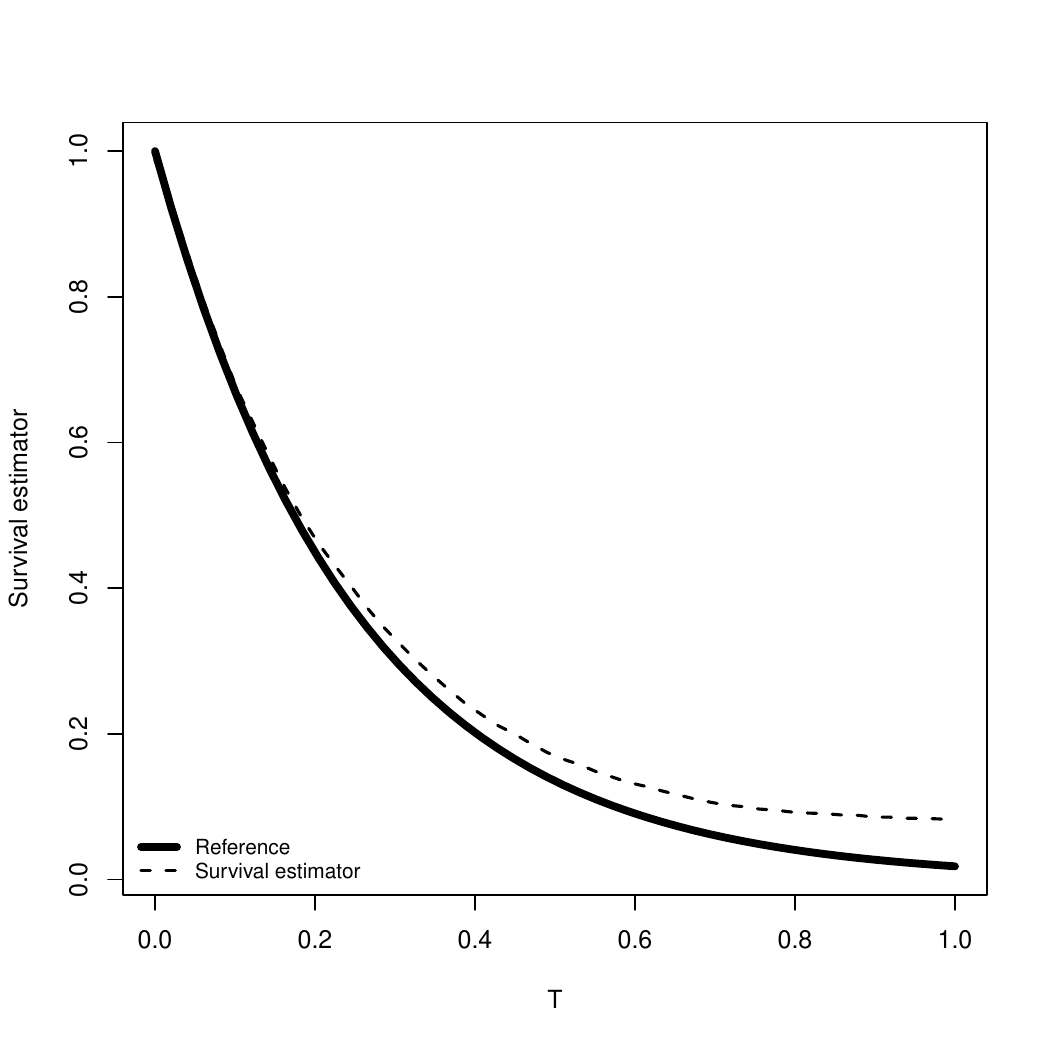}
\includegraphics[width =6cm, angle=0]{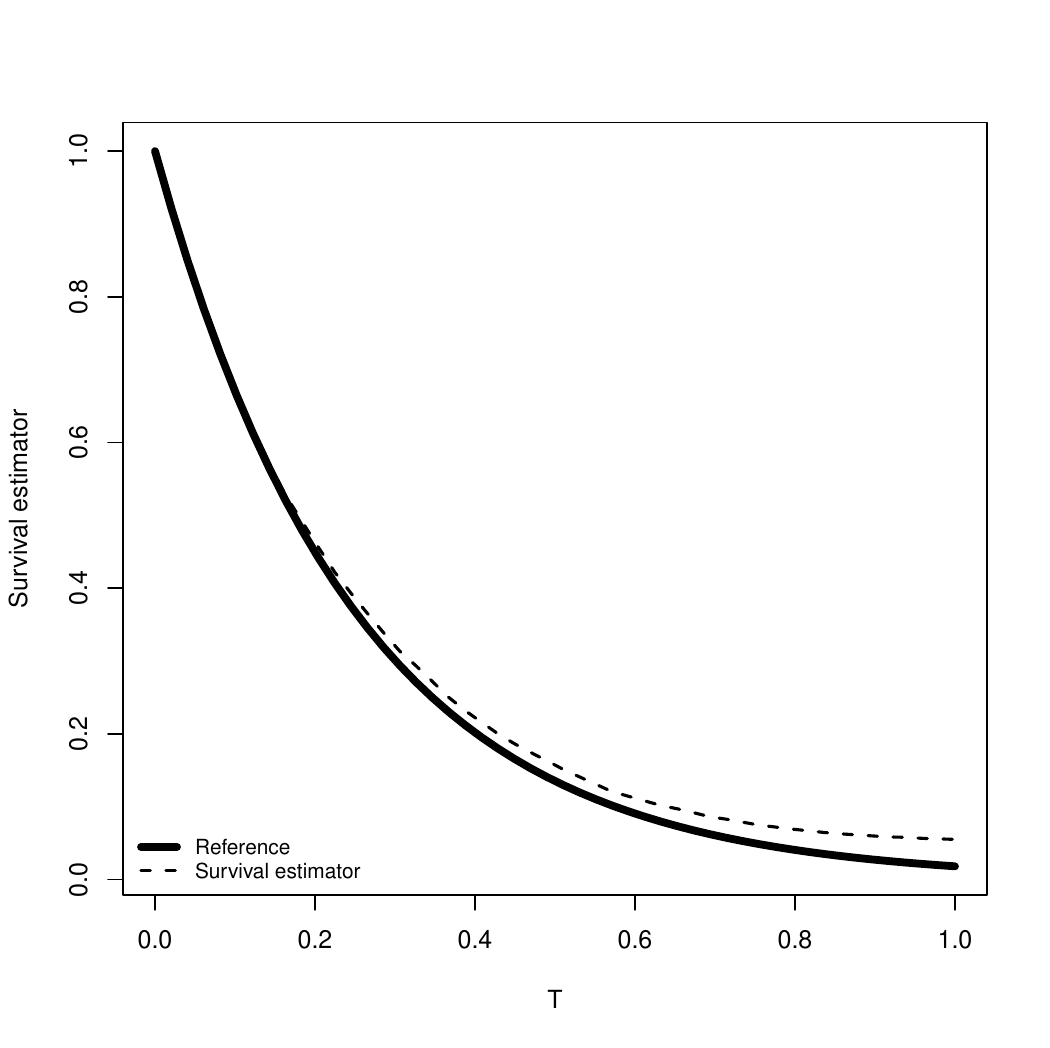}
\includegraphics[width =6cm, angle=0]{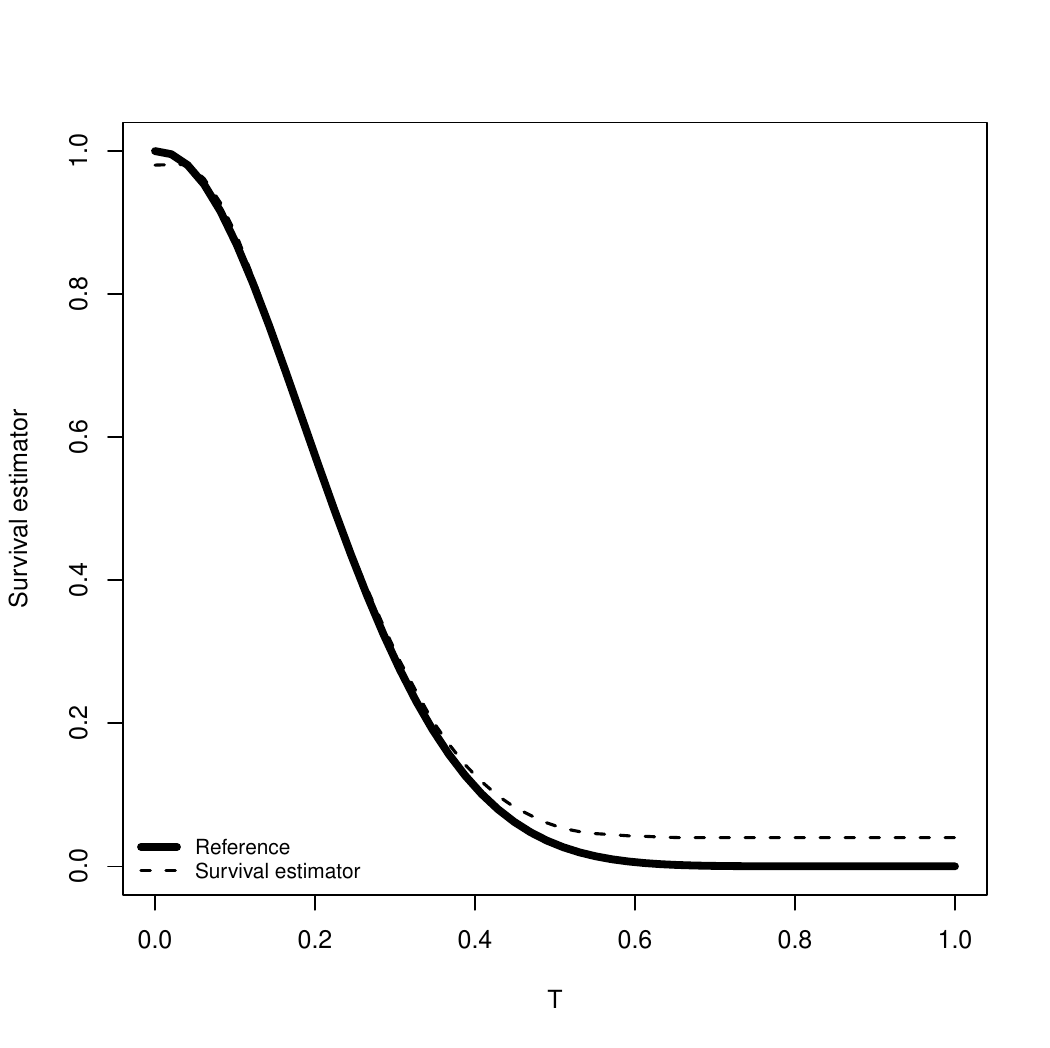}
\includegraphics[width =6cm, angle=0]{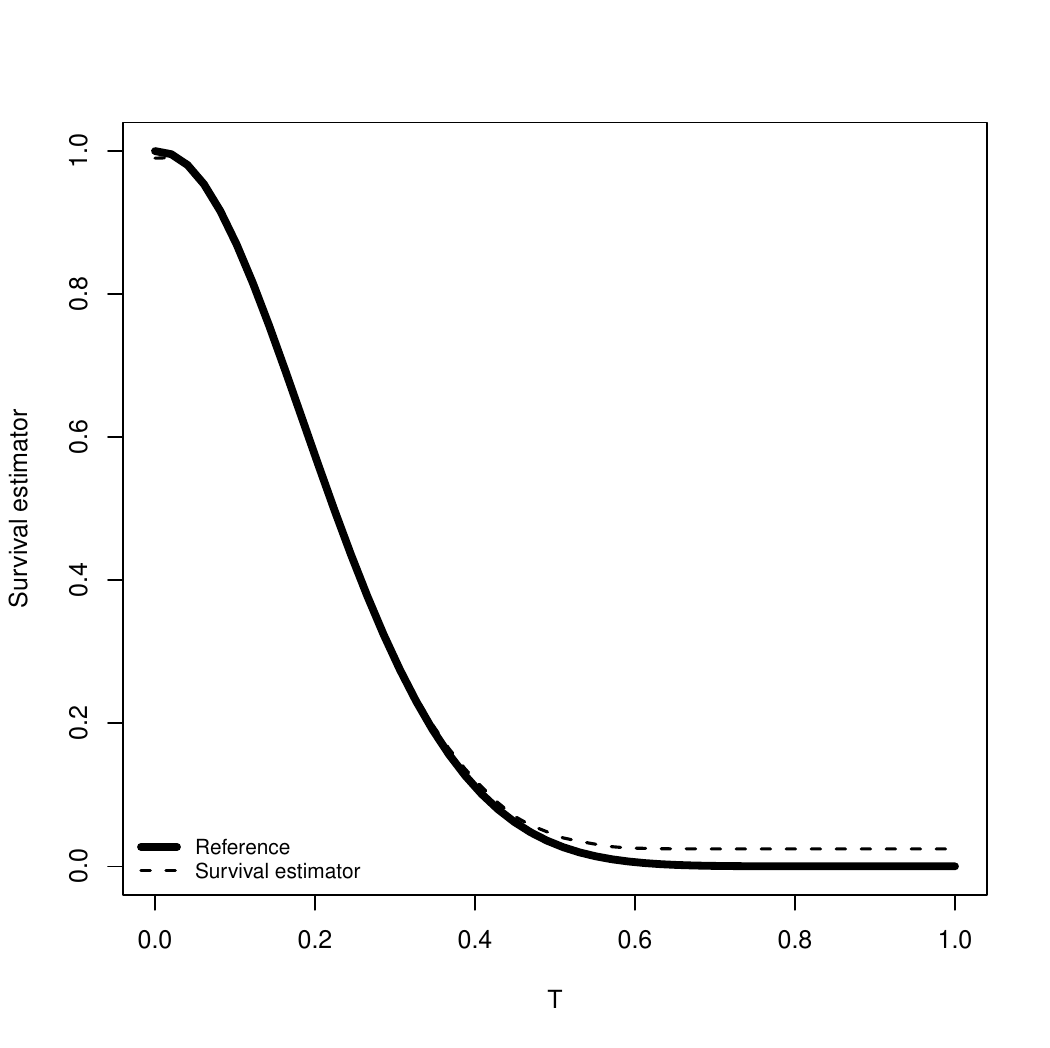}
\caption{Survival estimators $\overline{F}_{T,n}$, for $N=100$ samples  of size $n=100$ (left) and of size $n=200$ (right) for the models $a$ (up) and $b$ (down).}\label{Fig:marginals}
\end{figure}

\begin{figure}[H]
\centering
\includegraphics[width = 6cm,angle=0]{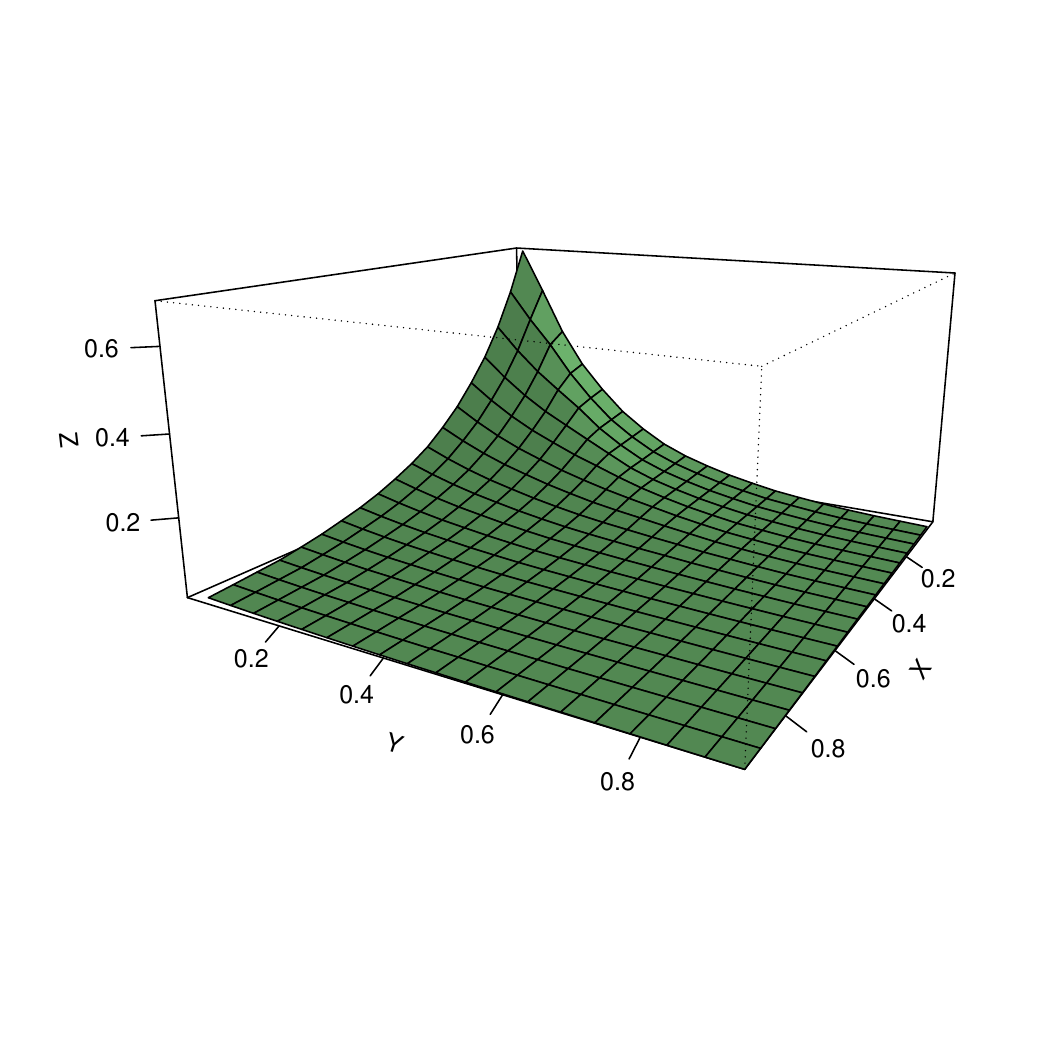}
\includegraphics[width = 6cm, angle=0]{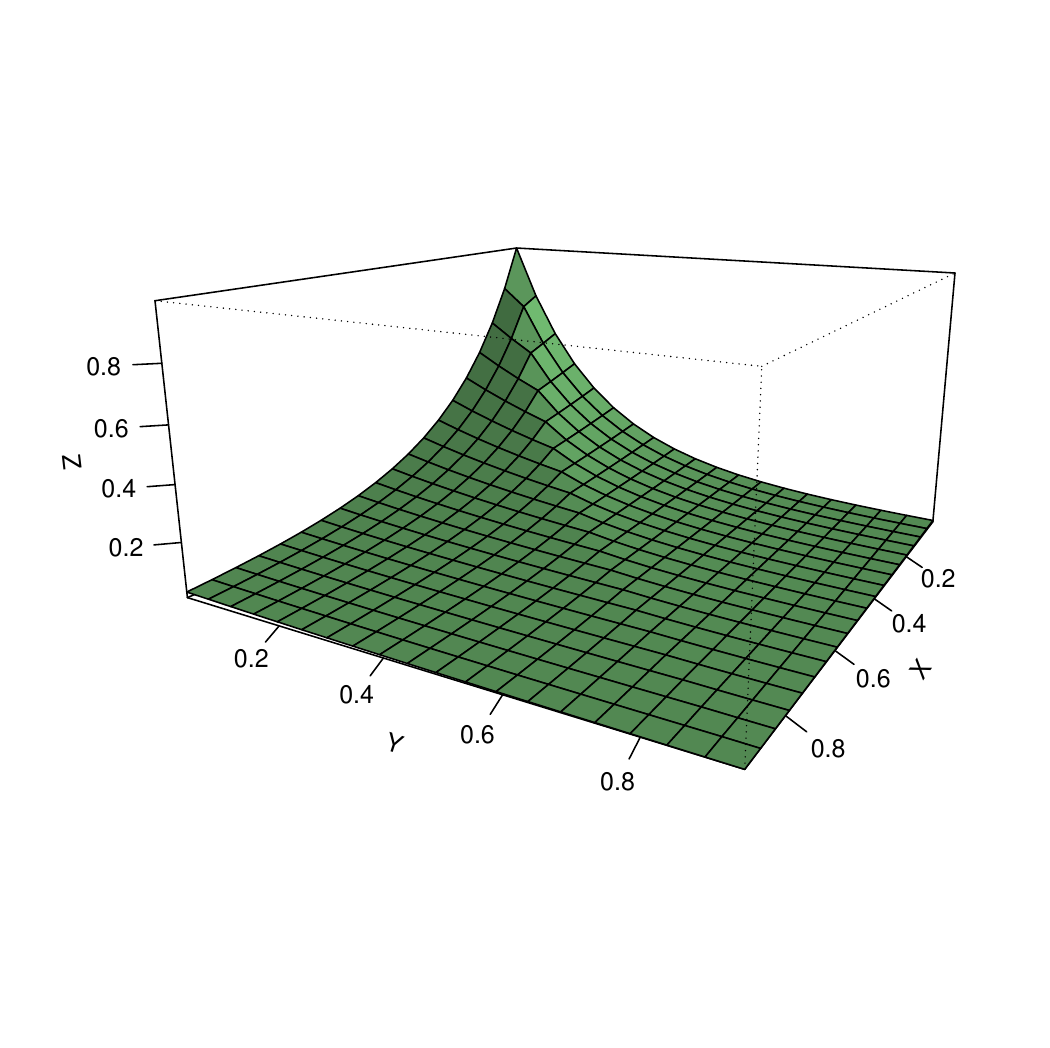}
\includegraphics[width = 6cm,angle=0]{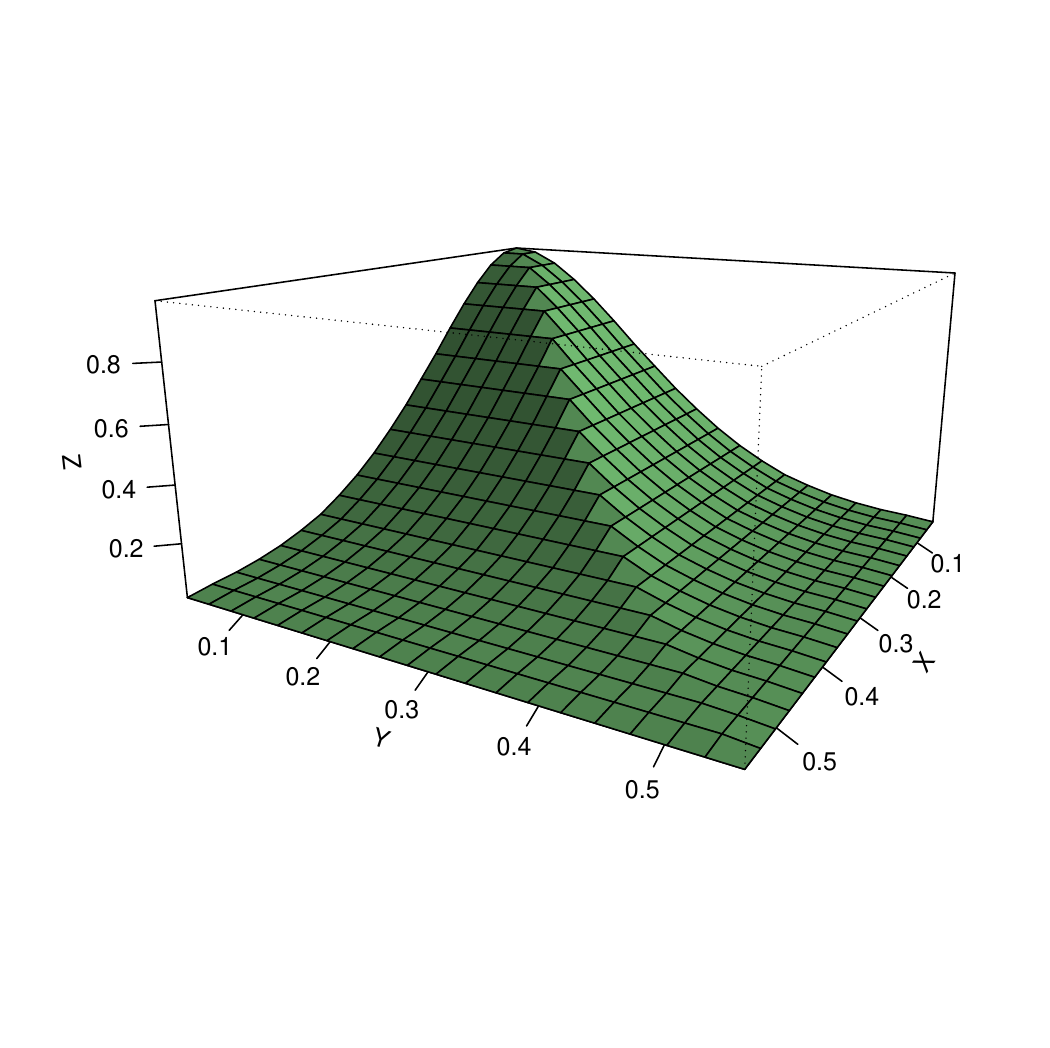}
\includegraphics[width = 6cm, angle=0]{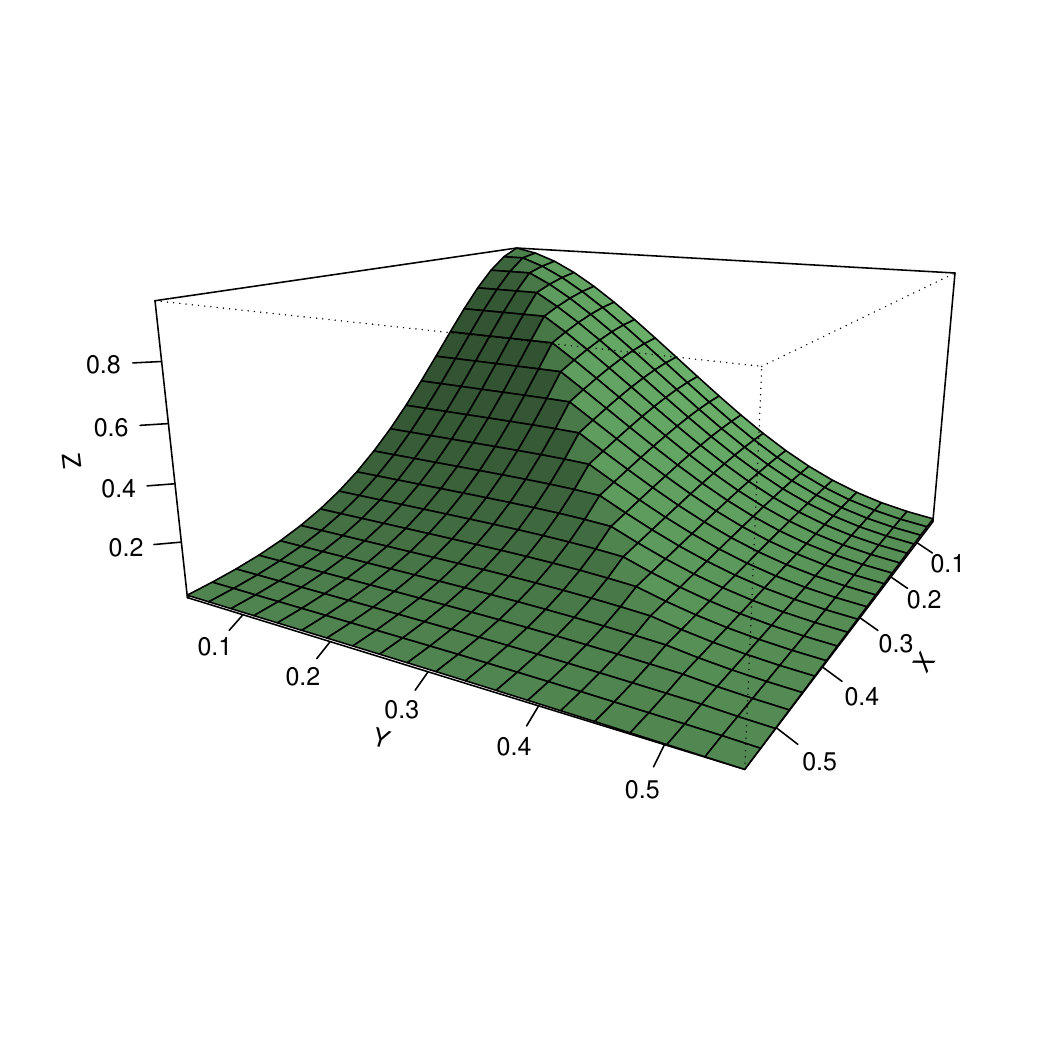}
\caption{Survival joint law estimator $\widetilde{P}_n$, for $N=100$ samples  of size $n=200$ (left) and the true survival joint law (right) for the model $a$ (up) and $b$(down). }\label{Fig:3d}
\end{figure}

In this subsection, we  ran a small simulation experiment to evaluate the performance of the survival  distribution estimator $\overline{F}_{T,n}$, the survival joint law  estimator $\widetilde{P}_n$  and the Kendall's tau $\tau_n$. We considered $N=100$ samples of sizes $n=50,100,200$ under the following models: 
\begin{enumerate}
\item[$a$ -]  The exponential distributions: $X_1 \sim \mathcal{E}(1)$, $X_2 \sim \mathcal{E}(2)$ and $X_{3} \sim \mathcal{E}(3)$,
\item[$b$ -]  The  beta distributions: $X_1 \sim \mathcal{B}(2,3)$, $X_2 \sim \mathcal{B}(10,10)$ and $X_{3} \sim \mathcal{B}(2.5,6)$,
\end{enumerate} 
where $\mathcal{E}$ and  $\mathcal{B}$ denote the  exponential and  Beta  distributions respectively. Figure \ref{alpha}  displays the average results for the estimation of the functions $\alpha_1$ and $\alpha_2$ under  the two models $a$ and $b$.  Figure \ref{Fig:marginals}  shows the mean of the survival estimator curves  while Figure \ref{Fig:3d}  displays the average results of the graphs for the survival copula estimation. As a qualitative measure of the efficiency over the different models, we computed the average of the integrated squared error (ISE) and the Kullback Leibler divergence (KL) of the survival joint law estimator $\widetilde{P}_n$ given by 
$$\overline{\text{ISE}}=\frac {1} {N} \sum_{k = 1}^ N \text{ISE}\left(\widetilde{P}_n^{(k)},\widetilde{P}\right) \quad\text{ and  }\,\overline{\text{KL}}=\frac{1}{N}\sum_{k = 1}^ N \text{KL}\left(\widetilde{P}_n^{(k)}|\widetilde{P}\right),$$ 
where $\widetilde{P}_n^{(k)}$ is the estimator computed from the $k^{th}$ sample and 
$$\text{ISE}\left(\widetilde{P}_n,\widetilde{P}\right) = \int_{[m,M]^2}\left(\widetilde{P}_k(x,y)-\widetilde{P}(x,y) \right)^2 dxdy, \quad \text{KL}(\widetilde{P}_n|\widetilde{P})=\int_{[m,M]^2}\widetilde{P}_n(x,y) \log \left(\frac{\widetilde{P}_n(x,y)}{\widetilde{P}(x,y)} \right)dxdy,$$ 
with $m=\min\limits_{1 \leq i \leq n}(Y_i)$ and $M=\max\limits_{1 \leq i \leq n}(Y_i)$. To investigate the performance of the proposed Kendall's tau estimator $\tau_n$, we computed the Bias and Mean squared error (MSE) of $\tau_n$ based on the true value of $\tau$.  It is worth noting  that the approximated value of the integral (\ref{eq::tau}) defines $\tau$ is obtained using the function \textsc{integrate} and the component \textsc{value} from the software \texttt{R}. From Table \ref{Tableltausim}, it can be concluded that $\tau_n$ is close to the true value of $\tau$.
 Figure \ref{Fig:marginals} shows that the survival estimators well behave and converge toward the true marginal survival functions. This is again a likeable feature from the GMO model since the dependency between $T$ and $C$ does not affect the KM estimator in this scenario.  Figure \ref{alpha}  shows that the estimators of $\alpha_1$ and $\alpha_2$ are constants under the model $a$ for a wide range of $t$ values, which matches with the MO exponential distribution. Furthermore,  the same figure shows that theses estimators can also capture non-constant function estimates well as it is the case under  model $b$. As expected from the consistency in Proposition \ref{asympnorm}, the results displayed in Tables \ref{Table1} show that the values for $\overline{\text{ISE}}$ and $\overline{\text{KL}}$ decrease along with the sample size $n$. This particularly emphasizes that the discrepancy between $\widetilde{P}_n$ and $\widetilde{P}$ is reduced as the plots in Figure \ref{Fig:3d} suggest.  Overall,  it can be concluded that the proposed non-parametric estimators are robust and perform well under the GMO model as it is theoretically expected. 

\subsection{Real data  applications}\label{sub2}
In this section, we evaluated the practical performances of our methods with two applications on real datasets. The first experiment  investigated the risk of three competing events during several soccer games where we compared the results between our approach and a maximum likelihood estimation based on the exponential MO model. In the second experiment, we proposed an application on the employment time for people working in a court of justice.          

\subsubsection{An application to the UEFA Champion’s League}
\begin{figure}[H]
\centering
\includegraphics[width=6cm, angle=0]{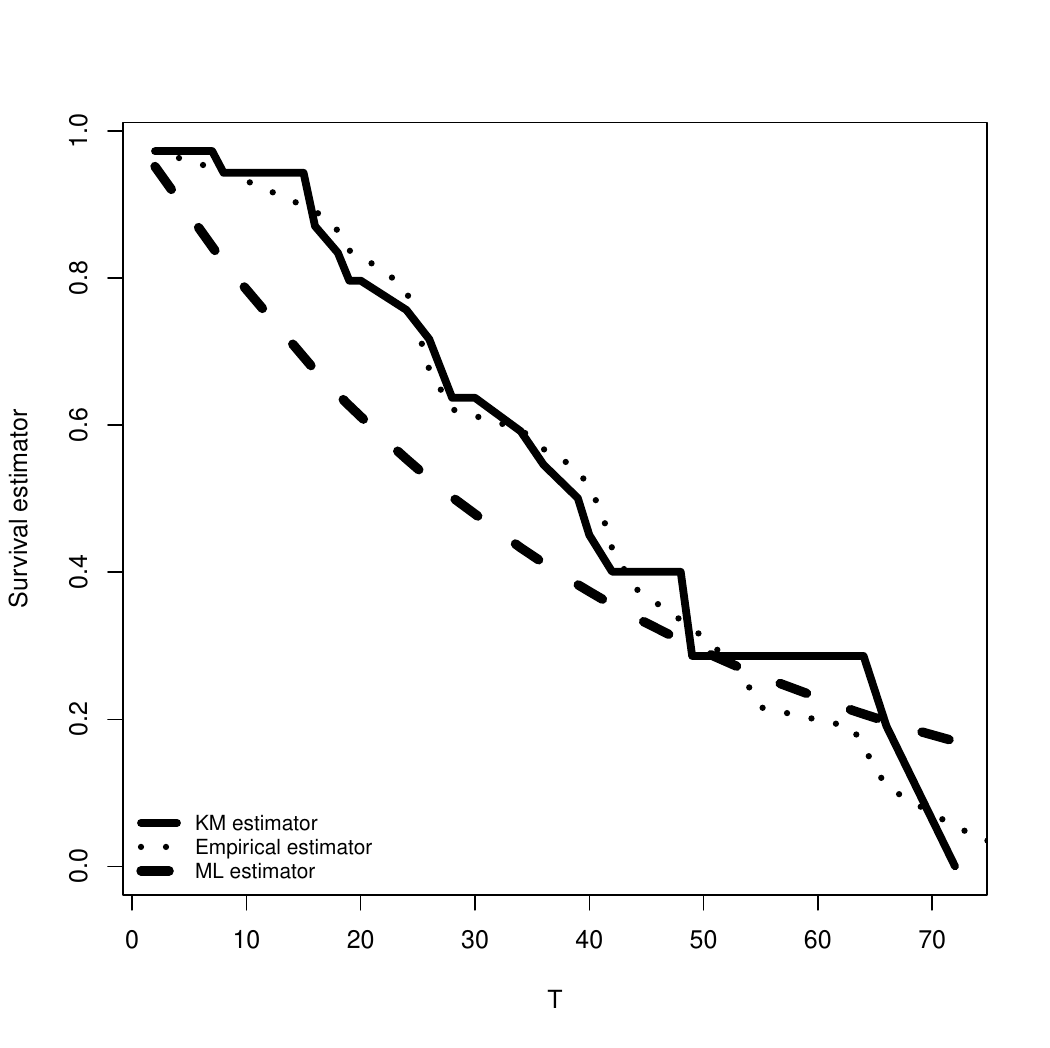}
\includegraphics[width=6cm, angle=0]{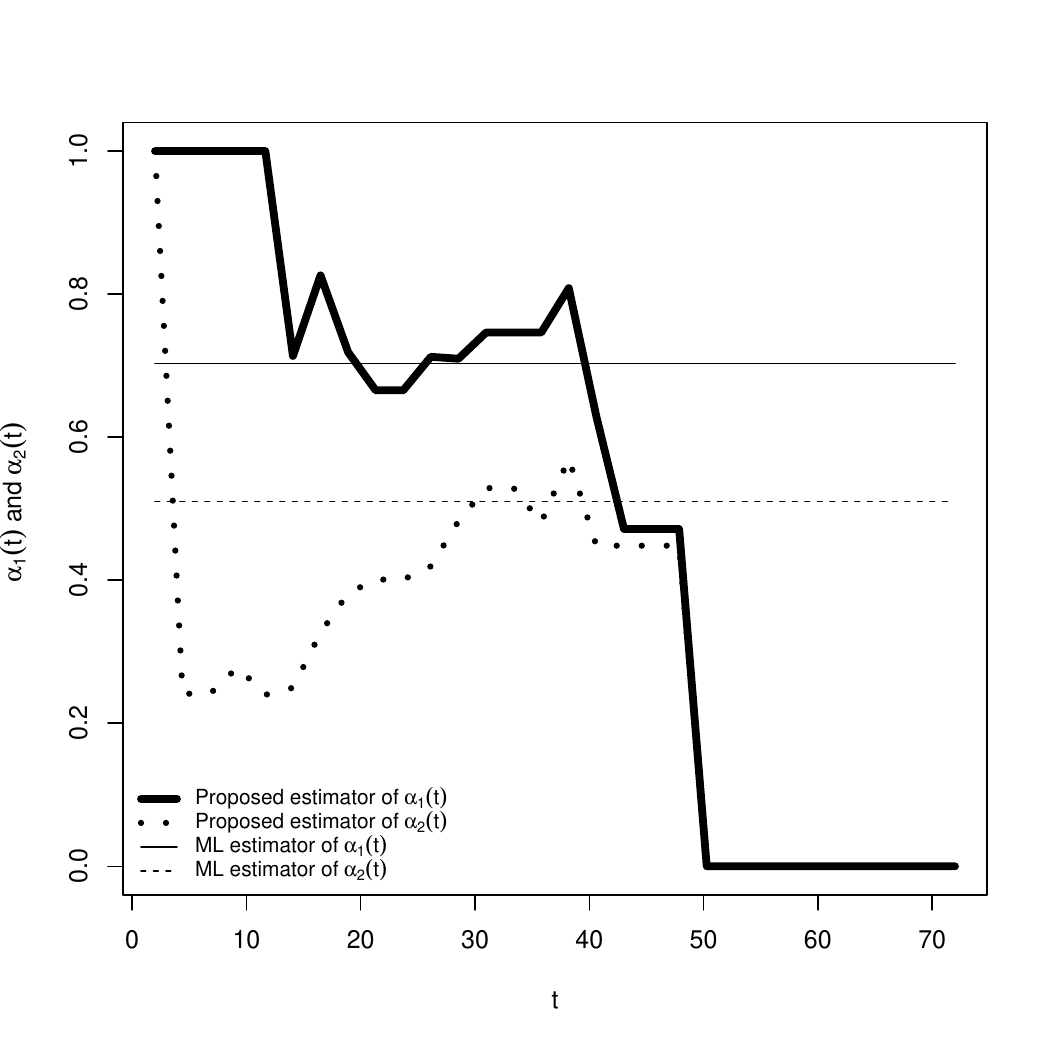}
\caption{Results for the survival distribution estimators (left) and the functions $\alpha_i$ (right) with the UEFA Champion’s League data.}\label{survivalfootanddatafoot}
\end{figure}

We  used soccer data from \cite{Mei07} which considers matches where there was at least one goal scored by the home team,  and at least one goal scored directly from a kick (penalty kick, foul kick or other) by any  of the  teams. Let $T$ be the time (in minutes) of the first goal scored by  one of the  teams, and $C$ be the time of the first goal of any type scored by the home team. At first,  the data makes it clear that we may have the events $T < C$, $T > C$ and $T = C$ for any of the  matches. According to \cite{Mei07}, these data may well have arisen from  an MO exponential distribution, implying that one shall consider three independent exponential random variables  $X_i\sim \mathcal{E}(\lambda_i)$ for $i=1,2,3$, with estimates given by $\widetilde{\lambda}_1=0.0073$, $\widetilde{\lambda}_2=0.0166$ and  $\widetilde{\lambda}_{3}=0.0173$. This yields the parametric survival marginal functions 
\begin{eqnarray*}
    \overline{F}_{T,ml}(t)=e^{-\left(\widetilde{\lambda}_1+\widetilde{\lambda}_3 \right)t}\quad\text{and}\quad\overline{F}_{C,ml}(t)=e^{-\left(\widetilde{\lambda}_2+\widetilde{\lambda}_3 \right)t},\quad t\geq 0,
\end{eqnarray*}
and
\begin{eqnarray*}
    \widetilde{P}_{ml}(t,s)=\min\left(\overline{F}_{T,ml}(s) \overline{F}_{C,ml}(t)^{1-\widetilde{\alpha_1}},\overline{F}_{T,ml}(t)\overline{F}_{C,ml}(s)^{1-\widetilde{\alpha_2}}\right),\quad t,s\geq 0,
\end{eqnarray*}
where $\alpha_i=\widetilde{\lambda}_{3}/(\widetilde{\lambda}_i+\widetilde{\lambda}_{3})$, $i=1,2$, for the joint parametric survival estimator. In this particular application, all the events  were surely observed at the end of the match, meaning that one can define the empirical survival estimators $\overline{F}_{T,ml}$ and $\overline{F}_{C,ml}$ together with the empirical survival joint distribution as follows:
\begin{eqnarray*}
    \widetilde{P}_{em}(t,s)=\frac{1}{n} \sum_{i=1}^n \mathds{1}_{\{ T_i > t, C_i > s\}},\quad t,s\geq 0.
\end{eqnarray*}
Likewise in section \ref{sec:num}, we proposed to compare the performances by computing the $\text{ISE}(\widehat{P},\widetilde{P}_{em})$ and $\text{KL}(\widehat{P}|\widetilde{P}_{em})$ scores for any considered joint law estimator $\widehat{P} \in \{ \widetilde{P}_n, \widetilde{P}_{ml}\}$. The returned values  were $\text{ISE}(\widetilde{P}_n,\widetilde{P}_{em})=0.001463$, $\text{ISE}(\widetilde{P}_{ml},\widetilde{P}_{em})=0.013737$ and $\text{KL}(\widetilde{P}_n|\widetilde{P}_{em})=0.022637$, $\text{KL}(\widetilde{P}_{ml}|\widetilde{P}_{em})=0.087960$. Our approach always shows the lowest values which supports our GMO model against the MO one. The same conclusion  was drawn in view of the distribution of $T$ in Figure \ref{survivalfootanddatafoot} where the ¨KM curve  was the closest to the empirical one. Moreover, the results  seem to contradict the idea that the data come from an actual MO model. In particular,   with the estimation of the functions $\alpha_i$, Figure \ref{survivalfootanddatafoot}  shows  that the results  are  not constant and do not always reach the expected values for the MO model of $\alpha_{i}^{ml}=\widetilde{\lambda}_3/\left(\widetilde{\lambda}_3+\widetilde{\lambda}_i\right)$. Overall, the discrepancy between our method and the empirical distribution tends to validate the GMO model, and rejecting,  at the same time,  the assumption of an exponential MO model since this latter  seems to rather a particular case of the GMO model.

\subsubsection {An application to the field of judges}
\begin{figure}[H]
\centering
\includegraphics[width=6cm, angle=0]{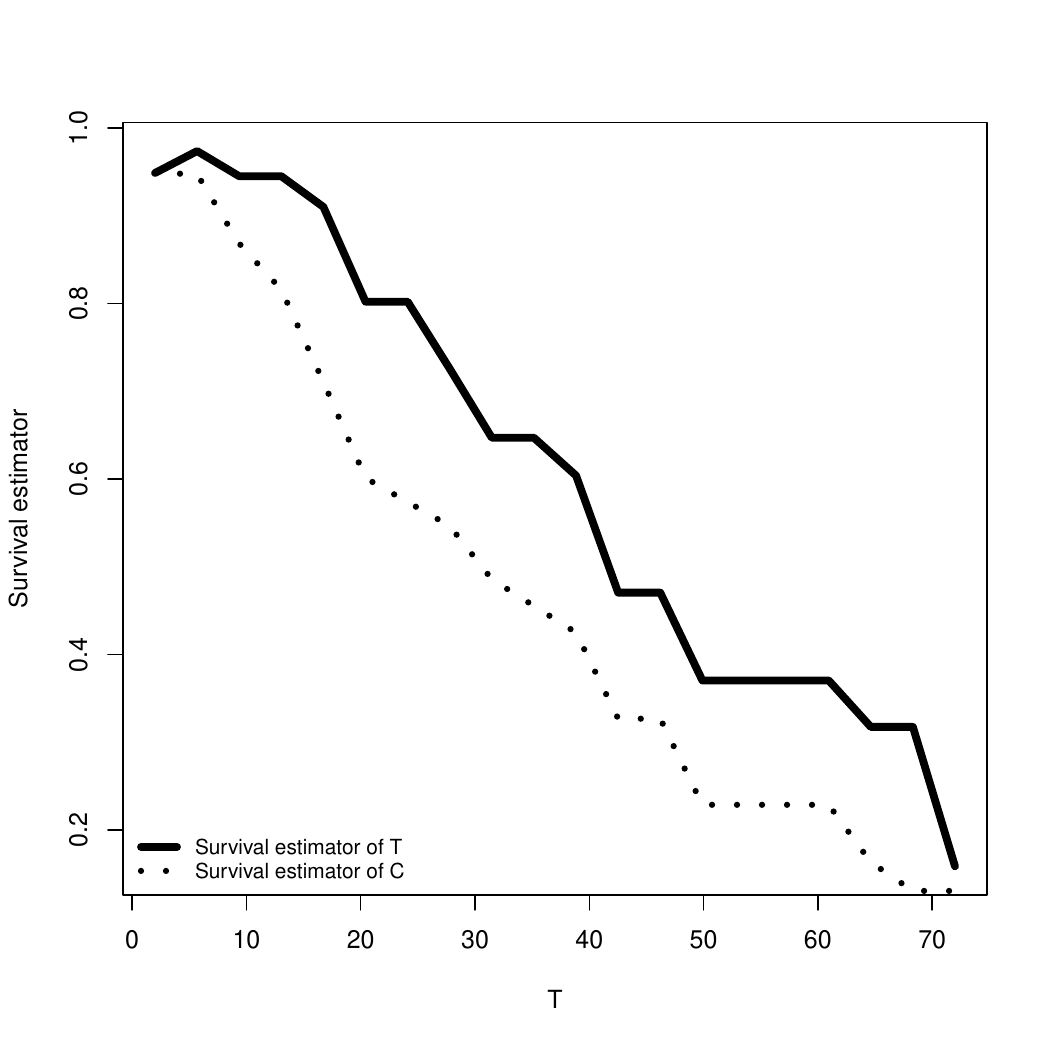}
\includegraphics[width=6cm, angle=0]{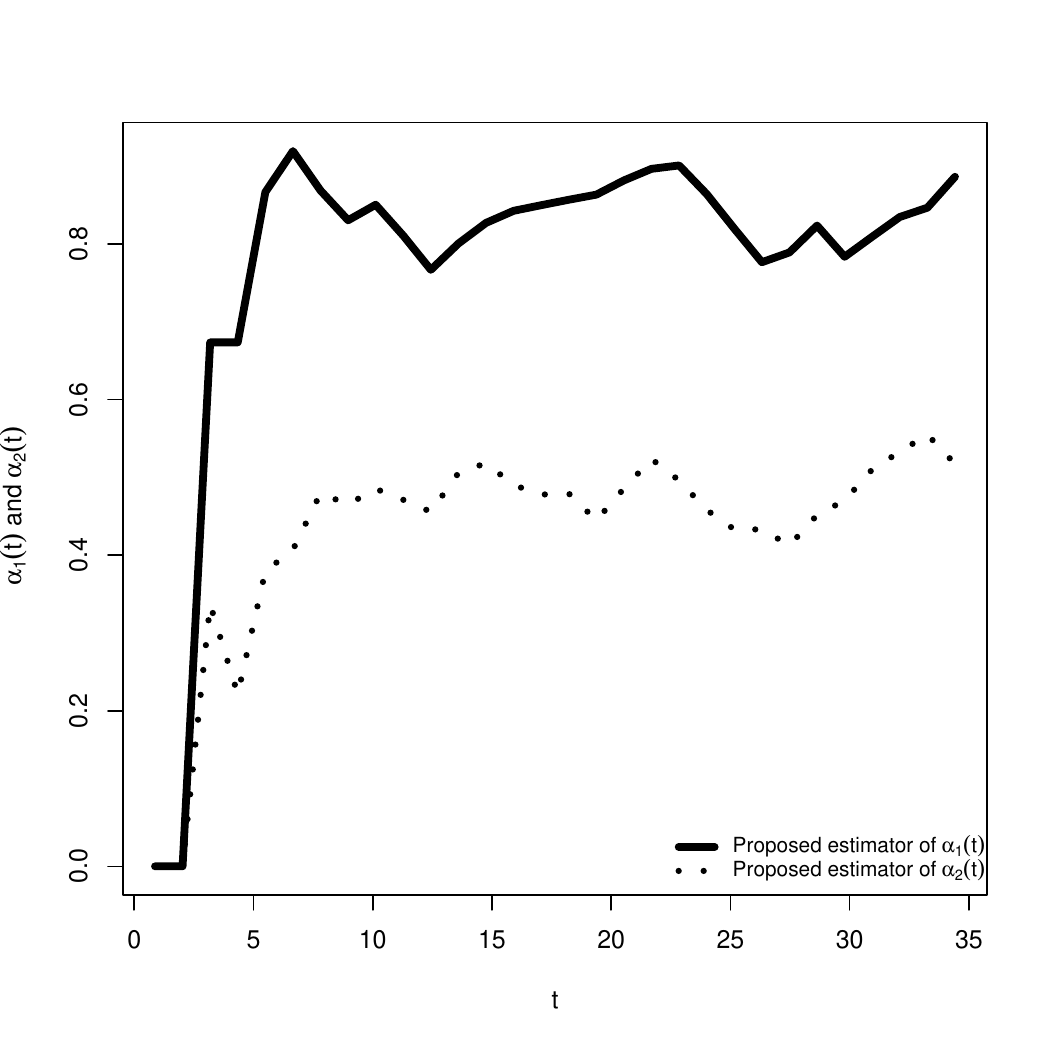}
\caption{Survival distribution estimators $\overline{F}_{T,n}$ and $\overline{F}_{C,n}$ for the judges data.}
\label{kmjustice}
\end{figure}

In this experiment, we considered an application to judges data\footnote{The data are available at \texttt{https://grodri.github.io/datasets/justices.csv}} which contain the length of service for workers of the U.S. Supreme Court  judges, treating death and retirement as competing risks. The sample consists of 113 individuals, and for each individual it records the time length of the tenure in the court service and assigns an indicator with the value 0 for the employees still in service, 1 for those who died while  on duty and 2 for those who resigned or retired. Here, we considered that the times of exercise until death $T$ and retirement $C$ are of interest, implying that one individual may express the event $T=C$ if she/he is still working at the record.  Figure \ref{kmjustice}  displays the results for the marginal survival estimators. The survival curves in Figure \ref{kmjustice}  imperatively suggest some proportionality assumption between the distributions of $T$ and $C$. This feature is also supported with the plots of the functions $\alpha_i$ which seem to remain constant for a large range of $t$ values. According to the Proposition \ref{prop::generalform}, we could thus infer that the cumulative hazard functions associated to the model are proportional. Intuitively, this could mean that the risks for an employee to quit the court  while still active, due to either a retirement plan or  death,  differ linearly. 

\section{Conclusion}\label{sec:con}
In this paper, we proposed a model and a non-parametric approach for survival times with dependent censoring based on the MO competing risk model. Several theoretical results  related to  the model function, in particular to the survival copula and the cumulative hazard functions,  were investigated. We  have also proposed estimators of the cumulative hazards and the survival joint law and shown their consistencies and asymptotic properties. The practical applicability of our method  was illustrated with the two studies on real datasets. We would also claim that this synthesis can be regarded as  a preliminary study  for further  investigations  of the survival copula based on the Marshall-Olkin's mechanism. Indeed, our work provides a theoretical  foundation of the characterization of the survival copula in the case of an extended Marshall Olkin model \citep[see][]{Pin15}, where only $(X_1,X_2)$ is independent of $X_3$. Another future potential research  involves the extension of our findings to the setting of the survival copula estimator using kernel approaches and Bernstein polynomials.

\section{Proofs}\label{sec:pro}

\textbf{Proof of proposition \ref{prop::generalform}} : In order to derive the survival copula, we start by expressing the joint distribution in terms of the functions $\widetilde{\Lambda}_i$ such that for any $s,t\geq 0$
\begin{eqnarray*}
 \mathbb{P}(T>t, C>s)&=& \mathbb{P}(X_1>t)\mathbb{P}(X_2>s)\mathbb{P}(X_{3}> t\vee s)
 \\&=& \exp(-\widetilde{\Lambda}_1(t)-
 \widetilde{\Lambda}_2(s)-\widetilde{\Lambda}_{3}(t\vee s)).
\end{eqnarray*}
With the relation $\max(x,y)=x+y–\min(x,y)$, we obtain that
\begin{eqnarray*}
\mathbb{P}\left(T>t,C>s \right)&=&\exp(-\widetilde{\Lambda}_1(t)-\widetilde{\Lambda}_2(s)-\widetilde{\Lambda}_{3}(t)-\widetilde{\Lambda}_{3}(s)+ \widetilde{\Lambda}_{3}(t\wedge s))
\\&=&\overline{F}_T(t)\overline{F}_C(s) \min\left(e^{\widetilde{\Lambda}_{3}(t)},e^{\widetilde{\Lambda}_{3}(s)}\right).
\end{eqnarray*}
Let $\alpha_i=\widetilde{\Lambda}_{3}/(\widetilde{\Lambda}_i+\widetilde{\Lambda}_{3})$ for i=1,2. Hence, we get for any $u\geq0$
\begin{eqnarray*}
\exp\left(\widetilde{\Lambda}_{3}(u)\right)&=&\overline{F}_T(u)^{-\alpha_1(u)}=\overline{F}_C(u)^{-\alpha_2(u)},
\end{eqnarray*}
which gives when  replaced in the joint distribution expression 
\begin{eqnarray*}
\tilde{C}(\overline{F}_T(t),\overline{F}_C(s))&=& \min \left(\overline{F}_C(s)\overline{F}_T(t)^{1-\alpha_1(t)}, \overline{F}_T(t) \overline{F}_C(s)^{1-\alpha_2(s)} \right).
\end{eqnarray*} 
We  can conclude the proof noting that for any generic distribution function $F$, $t=\overline{F}^{-}(\overline{F}(t))$ for any $t$ on the subset of $\mathbb{R}_+$ where $\overline{F}$ is continuous and strictly decreasing. \\

\noindent
\textbf{Proof of proposition \ref{prop::kendall}} : Let $\widehat{\Lambda}=\sum_{i=k}^3\Lambda_k$. Since the $\Lambda_i$ are striclty increasing, we know that they admit a first derivative on $\mathbb{R}_+$ almost everywhere. According to  Equation (30) in \cite{Mul18}, we have 
\begin{eqnarray*}
  \tau=4\int_0^{+\infty}e^{-2x}\Lambda_3\left(\widehat{\Lambda}^-(x)\right)dx.  
\end{eqnarray*}
By use of the change of variables $x=\widehat{\Lambda}(y)$ and the integration by parts formula, we  get
\begin{eqnarray*}
    \tau&=&4\int_0^{+\infty}e^{-2\widehat{\Lambda}(y)}[\widehat{\Lambda}(y)]'\Lambda_3(y)dy\\
    &=&4\left[-\dfrac{1}{2}e^{-2\widehat{\Lambda}(y)}\Lambda_3(y)\right]^{+\infty}_0+2\int_0^{+\infty}e^{-2\widehat{\Lambda}(y)}\Lambda_3'(y)dy\\
    &=&2\int_0^{+\infty}\overline{H}(y)^2d\Lambda_3(y)\\
    &=&2\int_0^{+\infty}\overline{H}(y)dH_3^1(y)
\end{eqnarray*}
where the last equality follows from  Proposition \ref{prop::kendall}.\\

\noindent
\textbf{Proof of Lemma \ref{lem::gmo}} : Let $E_1$, $E_2$ and $E_3$ be exponential random variables with parameters $\lambda_1=1/\alpha_1-1$, $\lambda_2=1/\alpha_2-1$ and 1 respectively. By definition of the MO model, this gives us that $(\min(E_1,E_3),\min(E_2,E_3))$ admit a survival copula $\widetilde{C}(u,v)=\min(u^{1-\alpha_1}v,uv^{1-\alpha_2})$, for $u,v\in[0,1]$. Since the copula definition for any random couple is free from the marginal distributions, it is  enough to independently transform $\min(E_1,E_3)$ and $\min(E_2,E_3)$ so that they have the correct distributions. Hence, we define
\begin{eqnarray*}
    T=F^-_X\left(1-e^{-(\lambda_1+1)\min(E_1, E_3)}\right)\quad\text{and}\quad C=F^-_Y\left(1-e^{-(\lambda_2+1)\min(E_2,E_3)}\right),
\end{eqnarray*}
which gives us
\begin{eqnarray*}
    X_1=F^-_X\left(1-e^{-(\lambda_1+1)E_1}\right),\quad X_2=F^-_Y\left(1-e^{-(\lambda_2+1)E_2}\right),\quad
    \text{and}\quad X_3=F^-_X\left(1-e^{-(\lambda_1+1)E_3}\right).
\end{eqnarray*}
It finally remains to show that $C$ follows the same law as $Y$, that is, we show that almost-surely
\begin{eqnarray*}
    F^-_Y\left(1-e^{-(\lambda_2+1)E_3}\right)=F^-_Y\left(1-\left[e^{-(\lambda_1+1)E_3}\right]^{\frac{\alpha_1}{\alpha_2}}\right)=F^-_X\left(1-e^{-(\lambda_1+1)E_3}\right),
\end{eqnarray*}
where the last equality comes from the proposition assumption since for any $t>0$
\begin{eqnarray*}
    F_Y^-(t)=\inf\{u>0:F_Y(u)\geq t\}&=&\inf\left\{u>0:1-e^{-\widetilde{\Lambda}_Y(u)}\geq t\right\}\\
    &=&\inf\left\{u>0:\widetilde{\Lambda}_Y(u)\geq -\log(
    1-t)\right\}\\
    &=&\inf\left\{u>0:\widetilde{\Lambda}_X(u)\geq -\dfrac{\alpha_2}{\alpha_1}\log(1-t)\right\}\\
    &=&F_X^-\left(1-(1-t)^{\frac{\alpha_2}{\alpha_1}}\right).
\end{eqnarray*}
\noindent
\textbf{Proof of proposition \ref{carac2}} : 
The proof is mainly based on results about the tail behavior of a distribution function in a maximum domain of attraction \citep[see for instance page 18 and Theorem 1.2.1 in][]{De06}. The idea is that $\overline{F}(t)$ roughly behaves like $t^{-1/\gamma}$  when $t\to+\infty$  and $\gamma>0$; likewise it roughly behaves like 
 $\overline{F}(x^{\star}-t)$ when  $t\downarrow 0$ and $\gamma<0$  with $x^{\star}$ its terminal point. Overall, we see that
\begin{eqnarray*}
    \widetilde{C}\left(u^{1/n},v^{1/n}\right)^n =\min\left(vu^{1-\alpha_1\left(\overline{F}^{-}_T(u^{1/n}) \right)},uv^{1-\alpha_2\left( \overline{F}^{-}_C(v^{1/n}) \right)}\right),\quad n\in\mathbb{N}^*.
\end{eqnarray*}
Note that  $\lim\limits_{n\to \infty} \alpha_1 \left( \overline{F}^{-}_T(u^{1/n}) \right)=x^{\star}_{F_T}$ and $\lim\limits_{n\to \infty} \alpha_2 \left( \overline{F}^{-}_C(u^{1/n}) \right)=x^{\star}_{F_C}$. In what follows, 
 the parameters for the limiting extreme value copula depend on the value of the limits $\lim_{t\uparrow x^{\star}_{F_T}}\alpha_1(t)$ and $\lim_{t\uparrow x^{\star}_{F_C}}\alpha_2(t)$. Hence, when $\gamma_i>0$ for $i=1,2$,  we have $x^{\star}_{F_T}=x^{\star}_{F_C}=\infty$ and 
\begin{eqnarray*}
\lim\limits_{t \to \infty}\dfrac{\widetilde{\Lambda}_i(t)}{\widetilde{\Lambda}_{3}(t)}&=& \lim\limits_{t \to \infty}\dfrac{\log(\overline{F}_i(t))}{\log(t)}\dfrac{\log(t)}{\log ( \overline{F}_{3}(t))}=\frac{\gamma_{3}}{\gamma_i},
\end{eqnarray*}

which implies that $\displaystyle{\lim\limits_{t \to \infty}\alpha_i(t)= \lim\limits_{t \to \infty}\frac{1}{1+\frac{\Lambda_i(t)}{\Lambda_3(t)}}= \frac{\gamma_i}{\gamma_i+ \gamma_3}}$. 
Thus, 
$$\lim\limits_{n \to \infty}  \widetilde{C}\left(u^{1/n},v^{1/n}\right)^n = \min \left( v u ^{1-\xi_1}, uv^{1-\xi_2}\right):=C_e,$$
where $\xi_i=\gamma_i/(\gamma_i+\gamma_3)$ for $i=1,2$. Hence, 
$\widetilde{C}$ belongs to the maximum domain of attraction of $C_e$. The case with $\gamma_i<0$ is similar.\\

\noindent
\textbf{Proof of proposition \ref{Lambda1etT}} : 
To prove the proposition, we develop the expression of $\displaystyle{\mathbb{P} \left(
Y \leq t, \delta^{(k)}=1 \right)}$. Since the arguments are similar for all $k=1,2,3$, we can prove the result for $k=1$ only. Clearly, we have  
\begin{eqnarray*}
H^1_1(t)&=& \mathbb{P}(X_1\leq t,X_1\leq X_2\leq X_{3})+ \mathbb{P}(X_1\leq t,X_1\leq X_{3}<X_2)
\\&=& \int_0^t \int_{[u,+\infty)}  \overline{F}_{3}(v^-) dF_2(v) dF_1(u) + \int_0^t \int_{[u,+\infty)}  \overline{F}_2(v) dF_{3}(v) dF_1(u)
\\&=& \int_0^t  \left[ \int_{[u,+\infty)}  \overline{F}_{3}(v^-) dF_2(v) + \int_{[u,+\infty)}  \overline{F}_2(v) dF_{3}(v) \right]  dF_1(u).
\end{eqnarray*}
Using the integration-by-part formula for the Stieltjes integral, we obtain for any $t\geq 0$
\begin{eqnarray*}
H^1_1(t)= \int_0^t \overline{F}_{3}(u^-) \overline{F}_2(u^-)dF_1(u)= \int_0^t \overline{F}_C(u^-)dF_1(u).
\end{eqnarray*}
This implies that $dH^1_1(u)/\overline{F}_C(u^-)=dF_1(u)$, which, when replaced in the definition for the cumulative hazard function, gives
\begin{eqnarray*}
\Lambda_1(t)=\int_0^t \frac{dF_1(u)}{\overline{F}_1(u^-)}=\int_0^t \frac{H^1_1(u)}{ \overline{H}(u^-)}.
\end{eqnarray*}
For $k=4$, since $\mathbb{P}(X_1=X_3)=0$, we note that $H_4^1=H_1^1+H_3^1$ and by use of the integration-by-part formula for the Stieltjes integral, we obtain
\begin{eqnarray*}
F_T(t)&=&\int_0^t\overline{F}_1(u^-)dF_3(u)+\int_0^t\overline{F}_3(u^-)dF_1(u)\\
&=&\int_0^t\frac{1}{\overline{F}_2(u^-)}dH_3^1(u)+\int_0^t\frac{1}{\overline{F}_2(u^-)}dH_1^1(u).
\end{eqnarray*}
Hence $dH_4^1(u)/\overline{F}_2(u^-)=dF_T(u)$ and it follows that
\begin{eqnarray*}
\Lambda_T(t)=\int_0^t \frac{dF_T(u)}{\overline{F}_T(u^-)}=\int_0^t \frac{H^1_4(u)}{\overline{H}(u^-)}.
\end{eqnarray*}
The proof is similar for $k=5$ and is thus omitted.\\ 

\noindent
\textbf{Proof of Proposition \ref{asympnormLambda}}
:  
Let $0 \leq t\leq \tau$. Based on Theorem 2 of \cite{Cai97}, we have the almost sure representation
\begin{eqnarray}\label{theo2cai}
\Lambda_{k,n}(t)=\Lambda_k(t)+\frac{1}{n}\sum_{i=1}^n\xi_{k,i}(t)+r_{k,n}(t)  \text{ for } k \in[\![1,5]\!],
\end{eqnarray}
where $\displaystyle{ \xi_{k,i}(t)=\sigma_{k,k}(Y_i,t) - \frac{\mathds{1}_{ \{ Y_i \leq t, \delta^{(k)}=1\}}}{\overline{H}(Y_i)} }$  and $\sup_{0\leq t\leq \tau}|r_{k,n}(t)|=o_\mathbb{P}(\sqrt{n})$. 
It is worth noting  that the work of \cite{Cai97} assumes that $T$ and $C$ are independent.
But  the proof of equation (\ref{theo2cai}) is valid in the case of independence or dependence between $T$
and $C$.  In fact, the result of Equation (\ref{theo2cai})  is  based on equation (7) in \cite{Cai97} which is expressed as follows: 
$$\Lambda_{k,n}(t)-\Lambda_k(t)=\int_0^t \left( \frac{1}{\overline{H}_n(s)}-\frac{1}{\overline{H}(s)}  \right) d \left[ H_{k,n}^1(s)-H_k^1(s) \right] \text{ for } k \in[\![1,5]\!].$$
 The rest of the proof of Equation (\ref{theo2cai}) is based on  the asymptotic expansion  of  the difference between the empirical estimators and their estimated 
$\displaystyle{\overline{H}_n(s)-\overline{H}(s) \text{ and } H_{k,n}^1(s)-H_k^1(s).}$
In view of Equation (\ref{theo2cai}), it  comes that the asymptotic behaviour of the random vector 
$$\left(\sqrt{n}(\Lambda_{k,n}(t_j)-\Lambda_k(t_j)),\,k\in[\![1,5]\!],j\in[\![1,m]\!]\right)$$
can be deduced from  that of $\left(n^{-1/2}\sum_{i=1}^n\xi_{k,i}(t_j),\,k\in[\![1,5]\!],j\in[\![1,m]\!]\right)$ which is a direct application of the multidimensional central limit theorem for $i.i.d.$ vectors. To  get the limit covariance matrix $\Sigma$, it then remains to compute the expectation of the products between the vector components. Note that almost surely 
\begin{eqnarray*}
\sigma_{k,k}(Y_1,t)&=&  \int_0^t \frac{\ind_{\{Y_1>u\}}}{\overline{H}(u)}d\Lambda_k(u).
\end{eqnarray*}
Let $0\leq t,s\leq \tau $ and $k,l\in[\![1,5]\!]$, we have
\begin{eqnarray*}
\mathbb{E}\left[\xi_{k,1}(t)\xi_{l,1}(s)\right]&=&\mathbb{E} \left[\sigma_{k,k}(Y,t) \sigma_{l,l}(Y,s)\right] -\mathbb{E} \left[ \sigma_{l,l}(Y,s)\frac{\mathds{1}_{ \{ Y \leq t, \delta^{(k)}=1 \}} }{\overline{H}(Y)}\right]\\
&&-\mathbb{E} \left[\sigma_{k,k}(Y,t)\frac{\ind_{\{Y\leq s,\delta^{(l)}=1 \}}}{\overline{H}(Y)} \right]+\mathbb{E}\left[\frac{\ind_{\{ Y \leq t, \delta^{(k)}=1\}}\ind_{\{ Y \leq s, \delta^{(l)}=1\}}}{\overline{H}(Y)^2} \right]\\
&=& \int_0^s \int_0^t  \frac{\overline{H}(u\vee v)}{\overline{H}(u)\overline{H}(v)}d\Lambda_k(u) d\Lambda_l(v)- \int_0^{s\wedge t} \frac{\Lambda_k(t)-\Lambda_k(u)}{\overline{H}(u)}d\Lambda_l(u)
\\ & & - \int_0^{s\wedge t} \frac{\Lambda_l(s)-\Lambda_l(u)}{\overline{H}(u)}d\Lambda_k(u)+ \mathbb{E}\left[\frac{\ind_{\{ Y \leq s\wedge t, \delta^{(k)}=\delta^{(l)}=1\}}}{\overline{H}(Y)^2} \right].
\end{eqnarray*}
When $k=l$, we observe that
\begin{eqnarray*}
\mathbb{E}\left[\xi_{k,1}(t)\xi_{k,1}(s)\right]=\int_0^s \int_0^t  \frac{\overline{H}(u\vee v)}{\overline{H}(u)\overline{H}(v)}d\Lambda_k(u) d\Lambda_k(v)- \int_0^{s\wedge t} \frac{\Lambda_k(t)+\Lambda_k(s)-2\Lambda_k(u)}{\overline{H}(u)}d\Lambda_k(u)+ \sigma_{k,k}(s,t)
\end{eqnarray*}
where
\begin{eqnarray*}
\int_0^s \int_0^t  \frac{\overline{H}(u\vee v)}{\overline{H}(u)\overline{H}(v)}d\Lambda_k(u) d\Lambda_k(v)=2\int_0^{s\wedge t}\frac{\Lambda_k(s\wedge t)-\Lambda_k(u)}{\overline{H}(u)}d\Lambda_k(u)+(\Lambda_k(s\vee t)-\Lambda_k(s\wedge t))\sigma_{k,k}(s,t)
\end{eqnarray*}
and for any $u>0$
\begin{eqnarray*}
\Lambda_k(t)+\Lambda_k(s)-2\Lambda_k(u)=2(\Lambda_k(s\wedge t)-\Lambda_k(u))+(\Lambda_k(s\vee t)-\Lambda_k(s\wedge t)).
\end{eqnarray*}
Hence
\begin{eqnarray*}
\int_0^s \int_0^t  \frac{\overline{H}(u\vee v)}{\overline{H}(u)\overline{H}(v)}d\Lambda_k(u) d\Lambda_k(v)=\int_0^{s\wedge t} \frac{\Lambda_k(t)+\Lambda_k(s)-2\Lambda_k(u)}{\overline{H}(u)}d\Lambda_k(u),
\end{eqnarray*}
and the proof is complete.\\

\noindent
\textbf{Proof of Proposition \ref{asympnorm}} : According to Lemma  1 of  \cite{Bre74}, the approximation of the KM estimator through the Nelson-Aalen estimator is negligible  and therefore either of them  can be used interchangeably. Formally, we have  
\begin{eqnarray*}
    0<-\log\left(\overline{F}_{k,n} (t)\right)- \Lambda_{k,n}(t)=O_\mathbb{P}(n^{-1}) \text{ for } k \in [\![1,5]\!],
\end{eqnarray*}
uniformly in $t\in[0,\tau]$ which means that  one can consider $\Lambda_{k,n}$ instead of $-\log\left(\overline{F}_{k,n}\right)$ with no changes on the asymptotic of $\widetilde{P}_n$. The result then is a direct consequence of the normality result from Proposition \ref{asympnormLambda} and the use of the Delta method when $\overline{F}_T(t)^{\alpha_1(t)}\neq\overline{F}_C(s)^{\alpha_2(s)}$.  For   $x,y,z\geq0$, we introduce the function $g(x,y,z)=\exp(-x)\exp(-y)^{1-z/y}$. For $n$ large enough, we observe that for any $t,s>0$
\begin{eqnarray*}
   && \hspace{-2cm}\sqrt{n}(\tilde{P}_n(t,s)-\tilde{P}(t,s))\\
    &=&
    \begin{cases}
       g(\Lambda_{4,n}(t),\Lambda_{5,n}(s),\Lambda_{3,n}(s))-g(\Lambda_T(t),\Lambda_C(s),\Lambda_3(s)) & \text{if } \overline{F}_T(t)^{\alpha_1(t)}<\overline{F}_C(s)^{\alpha_2(s)},\\[.2cm]
       g(\Lambda_{5,n}(s),\Lambda_{4,n}(t),\Lambda_{3,n}(t))-g(\Lambda_C(s),\Lambda_T(t),\Lambda_3(t)), & \text{if }\overline{F}_T(t)^{\alpha_1(t)}>\overline{F}_C(s)^{\alpha_2(s)},
    \end{cases}
\end{eqnarray*}
which weakly converges towards $\mathcal{N}(0,\Sigma(t,s))$ according to  the result of Proposition \ref{asympnormLambda} and  the Delta method applied to the function $g$, where $\Sigma(t,s)=Dg[\Lambda_T(t),\Lambda_C(s),\Lambda_3(s)] \Phi Dg[\Lambda_T(t),\Lambda_C(s),\Lambda_3(s)]^T$. Then, we get 
\begin{eqnarray*}
    \Sigma(t,s)=
    \begin{cases}
       \sigma_{3,3}(s,s)+\sigma_{4,4}(t,t)+\sigma_{5,5}(s,s)&
       \\+\,2(\sigma_{4,5}(t,s)-\sigma_{4,3}(t,s)-\sigma_{5,3}(s,s)), & \text{if } \overline{F}_T(t)^{\alpha_1(t)}<\overline{F}_C(s)^{\alpha_2(s)}.\\[.2cm]
       \sigma_{3,3}(t,t)+\sigma_{4,4}(t,t)+\sigma_{5,5}(s,s)&\\
       +\,2(\sigma_{4,5}(t,s)-\sigma_{4,3}(t,t)-\sigma_{5,3}(s,t)), & \text{if }\overline{F}_T(t)^{\alpha_1(t)}>\overline{F}_C(s)^{\alpha_2(s)}.
    \end{cases}
  \end{eqnarray*}
 In the case $\overline{F}_T(t)^{\alpha_1(t)}=\overline{F}_C(s)^{\alpha_2(s)}$, we consider the Skorokhod representation for the random vector $\left(\Lambda_{k,n}(u),\,k=3,4,5,u=s,t\right)$. Plainly, by the Taylor expansion we have
\begin{eqnarray*}
    \overline{F}_{4,n}(t)^{1-\alpha_{1,n}(t)} =\overline{F}_{T,n}(t)^{1-\alpha_{1,n}(t)}&=&\exp\left(-\Lambda_{4,n}(t) +\Lambda_{3,n}(t)\right)+O_\mathbb{P}(n^{-1})
    \\&=&\overline{F}_T(t)^{1-\alpha_1(t)}(1+n^{-1/2}[-N_4(t)+N_3(t)]+o_\mathbb{P}(n^{-1/2}))
\end{eqnarray*}
and
\begin{eqnarray*}
   \overline{F}_{5,n}(s)^{1-\alpha_{2,n}(s)} =\overline{F}_{C,n}(t)^{1-\alpha_{2,n}(t)}&=&\exp\left(-\Lambda_{5,n}(s) +\Lambda_{3,n}(s)\right)+O_\mathbb{P}(n^{-1})\\
    &=&\overline{F}_C(s)^{1-\alpha_2(s)}(1+n^{-1/2}[-N_5(s)+N_3(s)]+o_\mathbb{P}(n^{-1/2})),
\end{eqnarray*}
which yields
\begin{eqnarray*}
    \widetilde{P}_n(t,s)=\widetilde{P}(t,s)+n^{-1/2}\min\left(\overline{F}_T(t)^{1-\alpha_1(t)}[N_3(t)-N_4(t)],\overline{F}_C(s)^{1-\alpha_2(s)}[N_3(s)-N_5(s)]\right)+o_\mathbb{P}(n^{-1/2}),
\end{eqnarray*}
and the result follows.\\

\noindent
\textbf{Proof of Proposition \ref{prop::kendallest}} : 
We prove the asymptotic normality for $\tau_n$ by deriving its almost-sure representation. Seeing that based on (\ref{eq::tau}), we have
\begin{eqnarray*}
    \tau_n-\tau= 2\int_0^{+\infty}\overline{H}_n(x)dH_{3,n}^1(x)- 2\int_0^{+\infty}\overline{H}(x)dH_3^1(x),
\end{eqnarray*}
where 
\begin{eqnarray*}
    \overline{H}_n(x)=\dfrac{1}{n}\sum_{i=1}^n\ind_{\{Y_i>t\}}\quad\text{and}\quad H^1_{3,n}(x)=\dfrac{1}{n}\sum_{i=1}^n\ind_{\{Y_i \leq t,\delta^{(3)}_i=1\}}.
\end{eqnarray*}

By integration by part for the Stieltjes integral, this gives us
\begin{eqnarray*}
    \tau_n-\tau&=& 2\int_0^{+\infty}\overline{H}_n(x)-\overline{H}(x)dH_3^1(x)+2\int_0^{+\infty}H_{3,n}^{1}(x)-H_3^1(x)dH(x)+2\int_0^{+\infty}\overline{H}_n(x)-\overline{H}(x)\,d[H_{3,n}^1-H_3^1](x) 
    \\&=&\dfrac{2}{n}\sum_{i=1}^n\int_0^{+\infty}\ind_{\{Y_i>x\}}-\overline{H}(x)dH_3^1(x)+ 2\int_0^{+\infty}\ind_{\{Y_i \leq x,\delta^{(3)}_i=1\}}-H_3^1(x)dH(x) +R_n.
\end{eqnarray*}
We first have to prove that $R_n=o_\mathbb{P}(n^{-1/2})$ so that the asymptotic properties of the estimator $\tau_n$ are only carried by the left hand terms. We choose to use the strong representation results from Lemma 2 in \cite{Lo96}. We here only prove our result when $\inf\{t>0,H(t)=1\}=+\infty$, the proof for a finite support being similar and even simpler. We consider a positive non-random threshold $T_n$ and divide the interval $[0,T_n]$ into subintervals $[x_i,x_{i+1}]$ for $i=1,\ldots,k_n,$ where $k_n$ is positive sequence. Clearly, we have
\begin{eqnarray*}
    \dfrac{1}{2}|R_n|&\leq& \left|\int_0^{T_n}\overline{H}_n(x)-\overline{H}(x)\,d[H_{3,n}^1-H_3^1](x)\right|+ \left|\int_{T_n}^{+\infty}\overline{H}_n(x)-\overline{H}(x)\,d[H_{3,n}^1-H_3^1](x)\right|\\
    &=:& (R_{n,1}+R_{n,2}).
\end{eqnarray*}
We here select $T_n$ large enough so that the right-hand term becomes negligible. Indeed, on can see that
\begin{eqnarray*}
    R_{n,2}=\ind_{\{Y_{(n)}> T_n\}}R_{n,2}+\ind_{\{Y_{(n)} \leq T_n\}}\int_{T_n}^{+\infty}\overline{H}(x)\,d H_3^1(x),
\end{eqnarray*}
which implies that for $T_n=H^-({1-\alpha_n})$ (the $\alpha_n$-quantile of $H$) with $\alpha_n=o(n^{-1})$, we have $R_{n,2}\leq\int_{T_n}^{+\infty}\overline{H}(x)\,d H_3^1(x)\leq \alpha_n^2$ for $n$ large enough. We next consider the previous subdivision of the interval $[0,T_n]$ and choose $k_n=O(n^{1/2}l_n^{-1/2})$ where $l_n=\log(n)$, so that one can obtain 
\begin{eqnarray}
\label{eq::Rn1}
R_{n,1} &\leq &  2\max_{i=1,\ldots,k_n-1}\sup_{x\in[x_i,x_{i+1}]}|\overline{H}_n(x)-\overline{H}(x)-\overline{H}_n(x_i)+\overline{H}(x_i)|\\
\nonumber & & + k_n \sup_{0 \leq x \leq T_n} |\overline{H}_n(x)-\overline{H}(x)| \max _{1 \leq i \leq n-1} |H^1_{3,n}(x_{i+1}) -H^1_{3,n}(x_i)-H_3^1(x_{i+1}) +H_3^1(x_i) |.
\end{eqnarray}
We further divide each subinterval $[x_i,x_{i+1}]$ into smaller subintervals $[x_{i,j},x_{i,j+1}]$ for $j=1,\ldots,a_n$, where $a_n=O(n^{1/4}l_n^{-1/4})$, so that $|H(x_{i,j+1})-H(x_{i,j})|=O(T_nn^{-3/4}l_n^{3/4})=O(r_n)$ uniformly. It follows that for any $i=1,\ldots,k_n$
\begin{eqnarray*}
    && \hspace{-2cm}\sup_{x\in[x_i,x_{i+1}]}|\overline{H}_n(x)-\overline{H}(x)-\overline{H}_n(x_i)+\overline{H}(x_i)|\\
    &\leq&\max_{j=1,\ldots,a_n}\sup_{x\in[x_{i,j},x_{i,j+1}]}|\overline{H}_n(x)-\overline{H}(x)-\overline{H}_n(x_i)+\overline{H}(x_i)|\\
    &\leq& \max_{j=1,\ldots,a_n}\sup_{x\in[x_{i,j},x_{i,j+1}]}|\overline{H}_n(x)-\overline{H}(x_{i,j})-\overline{H}_n(x_i)+\overline{H}(x_i)|+O(r_n)\\
    &\leq&\max_{j=1,\ldots,a_n}|\overline{H}_n(x_{i,j})-\overline{H}(x_{i,j})-\overline{H}_n(x_i)+\overline{H}(x_i)|+O(r_n).
\end{eqnarray*}
According to the inequality in Lemma 1 from \cite{Lo96}
\begin{eqnarray*}
    \mathbb{P}\left(|\overline{H}_n(x_{i,j})-\overline{H}(x_{i,j})-\overline{H}_n(x_i)+\overline{H}(x_i)|>3r_n\right)\leq 2e^{-z}
\end{eqnarray*}
for $z\leq nr_n$ and $z\sigma^2\leq nr_n^2$ where $\sigma^2=O(T_n/k_n)$. Thus for $n$ large enough, one can choose $z=\dfrac{3}{4}l_n$ so that we have
\begin{eqnarray*}
    \mathbb{P}\left(\max_{i=1,\ldots,k_n-1}\max_{j=1,\ldots,a_n}|\overline{H}_n(x_{i+1})-\overline{H}(x_{i+1})-\overline{H}_n(x_i)+\overline{H}(x_i)|>3r_n\right)&\leq& 2k_na_ne^{-z}=l_n^{-3/4}
\end{eqnarray*}
where the right-hand term tends to zero. Since by the law of the iterated logarithm for the empirical process we have $\sup_{x>0}|\overline{H}_n(x)-\overline{H}(x)|=O(n^{-1/2}l_n^{1/2})$ almost surely \citep[see for instance][]{Vaart98}, the proof of the upper bound for the second component in (\ref{eq::Rn1}) follows the same arguments noting that $|H^1_3(x)-H^1_3(y)|\leq|H(x)-H(y)|$ for any $x,y>0$. This ensures that $|R_n|=O_\mathbb{P}(r_n)$ and means that it remains to prove that $r_n=o(n^{-1/2})$. 
When $\gamma=0$, the tail behavior of $H$ admits a exponential decay in the sense that for any $\beta>0$, $\lim_{t\to+\infty}t^\beta\overline{H}(t)=0$ \citep[see for instance p.18 in][]{De06}, which implies that $T_n=o(\alpha_n^{-1/\beta})$. Hence, for $\alpha_n=n^{-9/8}$ and $\beta=9$, we obtain
\begin{eqnarray*}
    r_n=T_n n^{-1/8}l_n^{3/4}n^{-5/8}=T_n\alpha_n^{1/\beta}l_n^{3/4}n^{-5/8}=o(n^{-1/2}).
\end{eqnarray*}
When $\gamma>0$, $\overline{H}$ is regularly varying with index $-1/\gamma$ which implies that $\alpha_n=\overline{H}(T_n)=T_n^{-1/\gamma}\ell(T_n)$ and $T_n=O(\alpha_n^{-\gamma})$. Define then $\delta=1/4-\gamma>0$ and choose $\alpha_n=n^{-3/2}$, then one can see that
\begin{eqnarray*}
    r_n=T_nn^{\gamma}l_n^{3/4}n^{-1/2-\delta}=n^{-\gamma/2}l_n^{3/4}n^{-1/2-\delta}O(1)=o(n^{-1/2}).
\end{eqnarray*}
It finally remains to derive the asymptotic variance of $\tau_n-\tau$. Let us define $\eta_i=2\int_0^{+\infty}\ind_{\{Y_i>x\}}-\overline{H}(x)dH_3^1(x)+ 2 \int_0^{+\infty}\ind_{\{Y_i>x,\delta^{(3)}=1\}}-H_3^1(x)dH(x)$ for $i=1\ldots n$. Since $(\eta_i)_{1 \leq i \leq n}$  are i.i.d random variables, it follows that 
$\text{Var}\left(\sqrt{n}(\tau_n-\tau)\right)=\mathbb{E}\left( \eta_1^2 \right)+o(1)$. Moreover, we have 
\begin{eqnarray*}
\mathbb{E}\left[ \eta_1^2 \right]&=& \mathbb{E}\left[ \left( 
 2\int_0^{+\infty}\ind_{\{Y>x\}}-\overline{H}(x)dH_3^1(x)+ 2\int_0^{+\infty}\ind_{\{Y \leq x,\delta^{(3)}=1\}}-H_3^1(x)dH(x) \right)^2\right]
\\&=&  \mathbb{E}\left[ \left( 2\int_0^{+\infty}\ind_{\{Y>x\}}-\overline{H}(x)dH_3^1(x) \right)^2\right]+ \mathbb{E} \left[ \left( 2\int_0^{+\infty}\ind_{\{Y \leq x,\delta^{(3)}=1\}}-H_3^1(x)dH(x)\right)^2 \right]
\\& &+8 \mathbb{E} \left[ \int_0^{+\infty}\ind_{\{Y>x\}}-\overline{H}(x)dH_3^1(x)  \int_0^{+\infty}\ind_{\{Y \leq x,\delta^{(3)}=1\}}-H_3^1(x)dH(x) \right]
\\&=:& I_1 +I_2 +I_3
\end{eqnarray*}
where
\begin{eqnarray*}
I_1 &=&4 \mathbb{E} \left[\int_0^{+\infty}  \int_0^{+\infty} \left(\mathds{1}_{\{Y>x\}}-\overline{H}(x) \right) \left(\mathds{1}_{\{Y>z\}}-\overline{H}(z) \right)  dH_3^1(x) dH_3^1(z)\right]
\\&=&  4\int_0^{+\infty}\int_0^{+\infty} \overline{H}(x\lor z) - \overline{H}(x) \overline{H}(z) dH_3^1(x) dH_3^1(z),\\
I_2 &=& 4 \mathbb{E} \left[ \int_0^{+\infty} \int_0^{+\infty} \left(\mathds{1}_{\{Y \leq x, \delta^{(3)}=1\}}-H_3^1(x) \right) \left(\mathds{1}_{\{Y \leq z, \delta^{(3)}=1 \}}-H_3^1(z) \right) dH(x) dH(z) \right]
\\&=&4\int_0^{+\infty}\int_0^{+\infty} H_3^1(x \wedge z) - H_3^1(x) H_3^1(z) dH(x) dH(z),\\
I_3&=&8 \mathbb{E} \left[ \int_0^{+\infty} \int_0^{+\infty} \left(\mathds{1}_{\{Y>x\}}-\overline{H}(x) \right)  \left(\mathds{1}_{\{Y \leq z, \delta^{(3)}=1\}}-H_3^1(z) \right) dH_3^1(x) dH(z)\right]
\\&=&8\int_0^{+\infty} \int_0^{+ \infty} \mathbb{P} \left(x < Y \leq z,\delta^{(3)}=1 \right)-\overline{H}(x) H_3^1(z)  dH_3^1(x) dH(z),
\end{eqnarray*}
and the result follows.

\section{Figures and table}

\begin{table}[h]
\centering
\begin{tabular}{c c c c c}
    \hline
    &  & $n$\\
    \hline
    & Model & 50 & 100 & 200\\
    \hline
    & $a$ & 0.008230 & 0.000657 & 0.000453  \\
    $\overline{\text{ISE}}$ &  $b$& 0.001592  & 0.001358 & 0.001206 \\
  \hline
& $a$ & 0.054817 & 0.008143 & 0.007572 \\
$\overline{\text{KL}}$ &  $b$ & 0.027034 & 0.026419 & 0.021504 \\
\hline
\end{tabular}
\caption{Results for $\overline{\text{ISE}}$ and $\overline{\text{KL}}$ for $N = 100$ trials of the proposed survival joint law estimator $\widetilde{P}_n$ for $n=50$, $n=100$ and $n=200$. }
\label{Table1}
\end{table}

\begin{table}[h]
\centering
\begin{tabular}{c c c c c}
    \hline
    &  & $n$\\
    \hline
    & Model & 50 & 100 & 200\\
    \hline
    $\tau$  & $a$ & $0.5000$  & $0.5000$  &  $0.5000$  \\
    &  $b$&  $0.6230$ & $0.6230$ & $0.6230$ \\
    \hline
    $\tau_n$  & $a$ & $0.5141$  & $0.4961$   &  $0.5009$ \\
    &  $b$& $0.6355$  & $0.6357$  & $0.6263$ \\
    \hline
 Bias  & $a$ & $0.0141$  & $0.0038$  & $0.0009$  \\
    &  $b$& $0.0125$  & $0.0127$ & $0.0033$ \\
  \hline
MSE  & $a$ & $0.0001$ &  $1.44e^{-5}$ & $8.78e^{-7}$  \\
  &  $b$ & $0.0001$  & $0.0001$ & $1.12e^{-5}$ \\
\hline
\end{tabular}
\caption{Results for $\text{Bias}$ and $\text{MSE}$ for $N = 100$ trials of the proposed Kendall's tau estimator $\tau_n$ for $n=50$, $n=100$ and $n=200$. }
\label{Tableltausim}
\end{table}

\begin{table}[h]
\centering
\begin{tabular}{c c c c}
\hline
&  UEFA data  & Justice data \\
 \hline
$\tau_n$   & $0.3433$  &   $0.4341$    \\
\hline
\end{tabular}
\caption{Results for $\text{Bias}$ and $\text{MSE}$ of the proposed Kendall's tau estimator $\tau_n$ for the real data.}
\label{Tabletaudata}
\end{table}



\end{document}